\input epsf

\documentclass[12pt,reqno]{amsart}
\usepackage{amscd}
\usepackage{color}
\usepackage{cite}
\usepackage{amssymb}

\input{epsf.sty}
\catcode `\@=11
\def\numberbysection{\@addtoreset{equation}{section}
         \renewcommand{\theequation}{\thesection.\arabic{equation}}}
\numberbysection
\def\subsubsection{\@startsection{subsubsection}{3}%
  \normalparindent{.5\linespacing\@plus.7\linespacing}{-.5em}%
  {\normalfont\bfseries}}
\newcommand\href[1]{\empty}

\def\cal{\mathcal}

\def\isoto{\stackrel{\sim}{\rightarrow}}
\def\a{\alpha}

\def\t{\tau}

\def\nn{\nonumber}
\def\la{\langle}
\def\ra{\rangle}

\def\Z{\mathbb Z}

\def\nn{\nonumber}

\def\text{\mbox}

\newtheorem{thm}{Theorem}[section]
\newtheorem{lem}[thm]{Lemma}
\newtheorem{prop}[thm]{Proposition}
\newtheorem{cor}[thm]{Corollary}

\theoremstyle{definition}

\def\CO{{\mathcal O_\mu}}

\def\Ch{Chain}
\def\A{\widetilde{Arc}}
\def\Rp{{\mathbb R}_{>0}}
\def\PB{\mathcal{P}({\mathcal B})}
\def\BV{\Delta}
\def\Ass{\mathcal A}

\def\pp{P^2_1}
\def\pa{P^A}
\def\punc{P^1_2}
\def\B{\mathcal B}

\begin{document}
\title{Closed/open string diagrammatics}

\author[Ralph M. Kaufmann and R.C. Penner]
{Ralph M. Kaufmann\\
Department of Mathematics\\
University of Connecticut,
Storrs, USA\\
~
\\
R. C.  Penner\\
Departments of Mathematics and Physics/Astronomy\\
University of Southern California,
Los Angeles, USA}

\email{kaufmann@math.uconn.edu \\
rpenner@usc.edu}

\begin{abstract}
We introduce a combinatorial model based on measured foliations in
surfaces which captures the phenomenology of open/closed string
interactions.  The predicted equations are
derived in this model, and new equations can be discovered as
well. In particular, several new equations together with known
transformations generate the combinatorial version of open/closed
duality.  On the topological and chain levels, the algebraic
structure discovered is new, but it specializes to a
 modular
bi-operad on the level of homology.\end{abstract}

\maketitle

\section*{Introduction}  There has been considerable activity
 towards a satisfactory
diagrammatics of open/closed string interaction and the underlying
topological field theories.
For closed strings on the topological level, there are the fundamental results of Atiyah and Dijkgraaf
\cite{Atiyah, Dij1}, which are nicely summarized
in \cite{Kock}.
The topological open/closed theory has proved to
be trickier
since there have been additional unexpected axioms, notably the Cardy
condition \cite{CL,Lew,Laz,Segal2,MS}; this algebraic
background is again nicely summarized in \cite{LP}.

In closed string field theory
\cite{Wittennoncomm,KKS,SaaZwi}, there are many new algebraic
features \cite{Zwiebach,LZ,LZ2}, in particular, coupling to gravity
\cite{topgravWi,phases} and a Batalin-Vilkovisky structure
\cite{Getz1,Getz2}. This BV structure has the same origin as that
underlying string topology \cite{CS,C1,C2,V,K} and the
decorated moduli spaces \cite{KLP,P2,woods}.

In terms of open/closed theories beyond the topological level,
many interesting results have been established for $D$-branes
\cite{Bru1,HIV,Douglas,DD,Cardy,KatzSharpe,orlov,Letal,Laz2,FS,KLi1,KLi2,FRS,Diac,HerbstLaz,Lazrec,KR,KLi3}
and Gepner models in particular \cite{GutSa,Bru2,Horec}.
Mathematically, there has also been work towards generalizing
known results to the open/closed setting
\cite{Costello,Harrelson,BCR,KS1,KS2}.

We present a model which accurately reflects the standard
phenomenology of interacting open/closed strings and which
satisfies and indeed rederives the ``expected'' equations
of open/closed topological field theory and the BV-structure of
the closed sector. Furthermore, the model allows the calculation
of many new equations, and there is an infinite algorithm for
generating all of the equations of this theory on the topological
level.  A finite set of equations, four of them new, are shown to
generate open/closed duality.

The rough idea is that as the strings move and interact, they form the leaves
of a foliation, the ``string foliation'', on their world-sheets.
Dual to this foliation is another  foliation of the world-sheets,
which comes equipped with the additional structure of a
``transverse measure''; as we shall see, varying the transverse
measure on the dual ``measured foliation'' changes the combinatorial type of
the string foliation.

The algebra of these string interactions  is then given by gluing
together the string foliations along the strings, and this
corresponds to an appropriate gluing operation on the dual
measured foliations. The algebraic structure discovered is
new, and we axiomatize it (in Appendix~A) as a ``closed/open'' or ``c/o structure''.
This structure is present on the topological level of string
interactions as well as on the chain level.  On the homology
level, it induces the structure of a modular bi-operad, which
governs c/o string algebras (see Appendix~A and Theorem~\ref{homthm}).

Roughly, a measured foliation in a surface $F$ is a collection of rectangles of some fixed widths and
unspecified lengths foliated by horizontal lines  (see Appendix B for the precise definition).  One glues such a
collection of rectangles together along their widths in the natural measure-preserving way (cf. Figure~5), so as to
produce a measured foliation of a closed subsurface of $F$. In the transverse direction, there is a natural foliation
of
each rectangle also by its vertical string foliation, but this foliation has no associated transverse measure.
In effect, the physical length of the string is the width of the corresponding rectangle.
A measured foliation does {\sl not} determine a metric on the surface, rather,
one impressionistically thinks of a measured foliation as describing {\sl half} of a metric
since the
widths of the rectangles are determined but not their lengths (see also \S B.1 for more details).

Nevertheless, there is a condition that we may impose on measured
foliations by rectangles, namely, a measured foliation of $F$ by
rectangles is said to {\it quasi fill} $F$ if every component of
$F$ complementary to the rectangles is either a polygon or an
exactly once-punctured polygon. The cell decomposition of
decorated Riemann's moduli space for punctured surfaces
\cite{HubMas,Strebel,Harstab,Harvir,P3} has been extended to
surfaces with boundary in \cite{P2}, and the space of quasi
filling measured foliations by rectangles again turns out to be
naturally homotopy equivalent to Riemann's moduli space of $F$
(i.e., classes of structures on surfaces with one distinguished
point in each hyperbolic geodesic boundary component; see the next
section for further details). Thus, in contrast to a measured
foliation impressionistically representing half a metric, a quasi
filling measured foliation actually {\sl does} determine a
conformal class of metrics on $F$. See the closing remarks for a
further discussion of this ``passage from topological to conformal
field theory''.

 \vskip .2in

\hskip .7in\epsffile{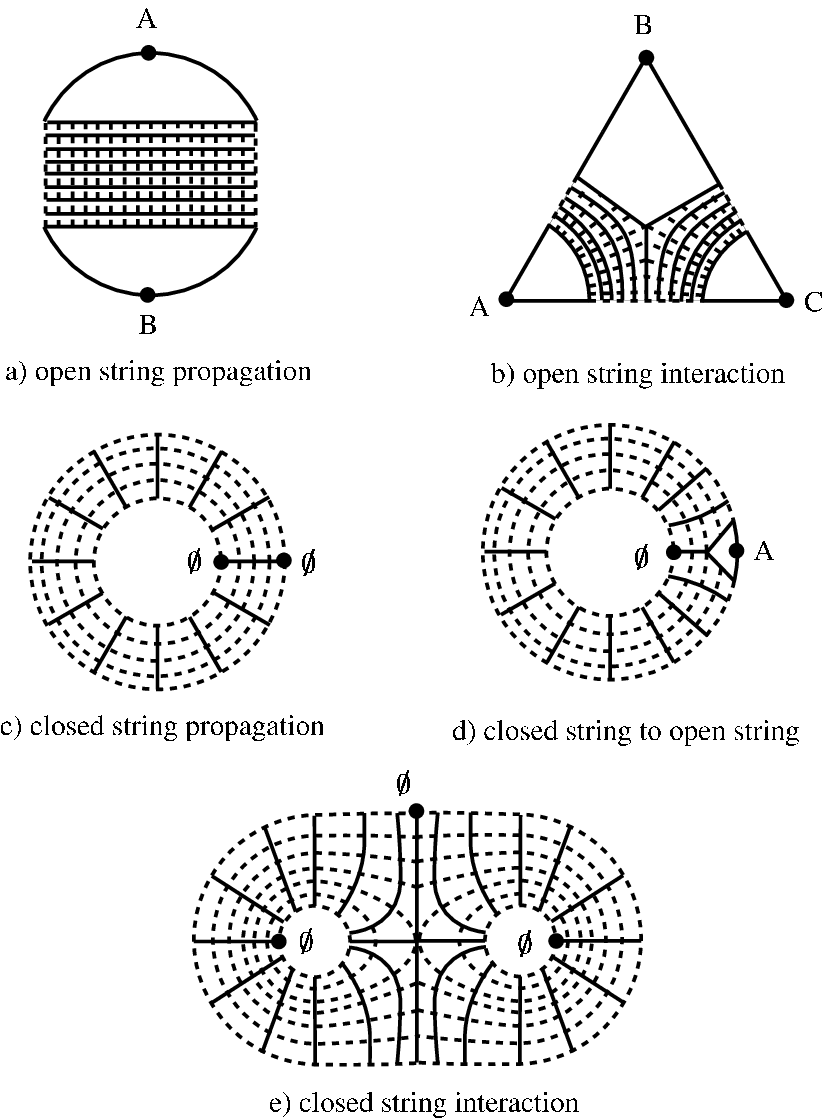}

\noindent \begin{footnotesize}{{{\bf Figure 1} Foliations for
several string interactions, where the strings are represented by
dashed lines and the dual measured foliation by solid lines. The
white regions in parts a-b are for illustration purposes only} }
\end{footnotesize}

\vskip .2in

More explicitly in Figure~1, each boundary component comes
equipped with a non-empty collection of distinguished points that
may represent the branes, and the labeling will be explained
presently.  That part of the boundary that is disjoint from the
foliation and from the distinguished points has no physical
significance: the physically meaningful picture arises by
replacing each distinguished point in the boundary by a small
distinguished arc (representing that part of the interaction that
occurs within the corresponding brane) and collapsing to a point
each component of the boundary disjoint from the foliation and
from the distinguished arcs.

Since the details for general measured foliations may obfuscate the relevant combinatorics and phenomenology of
strings,
we shall restrict attention for the most part to the special measured foliations where each non-singular leaf is an
arc properly embedded in the
surface.  The more general case is not without interest (see Appendix~B).

The natural equivalence classes of such measured foliations are in one-to-one correspondence with
``weighted arc families'', which are appropriate homotopy classes of properly and disjointly embedded arcs together
with
the assignment
of a positive real number to each component (see the next section for the precise definition).  Furthermore,
the quotient of this subspace of foliations by the mapping class group is closely
related to Riemann's moduli space of the surface (again see the next section for the precise statement).

A {\it windowed surface} $F=F^s_g(\delta_1,\ldots ,\delta _r)$ is a
smooth oriented surface of genus $g\geq 0$ with $s\geq 0$ punctures and
$r\geq 1$ boundary components together with the specification of a
non-empty finite subset $\delta _i$ of each boundary component,
for $i=1,\ldots ,r$, and we let $\delta =\delta
_1\cup\cdots\cup\delta _r$ denote the set of all distinguished points in the boundary
$\partial F$ of $F$ and let  $\sigma$ denote the set of all punctures.
The set of components of
$\partial F-\delta$ is called the set $W$ of {\it windows}.

In the physical context of interacting closed and open strings,
the open string endpoints are labeled by a set of branes in the
physical target, and we let ${\mathcal B}$ denote this set of
brane labels, where we assume $\emptyset\notin{\mathcal B}$.  In
order to account for all possible interactions, it is necessary to
label elements of $\delta\cup\sigma$ by the power set ${\mathcal
P}({\mathcal B})$ (comprised of all subsets of ${\mathcal B}$).  In effect,
the label $\emptyset$ denotes closed strings, and the label $\{
B_1,\ldots ,B_k\}\subseteq{\mathcal B}$ denotes the formal
intersection of the corresponding branes.  This intersection in
the target  may be empty in a given physical circumstance.

A {\it brane-labeling} on a windowed surface $F$ is a function
$$\beta :\delta\sqcup\sigma\to{\mathcal P}({\mathcal B}),$$
where $\sqcup$ denotes the disjoint union, so that if $\beta
(p)=\emptyset$ for some $p\in\delta$, then $p$ is the unique point
of $\delta$ in its component of $\partial F$.  A brane-labeling
may take the value $\emptyset$ at a puncture. (In effect,
revisiting windowed surfaces from \cite{Pthesis} now with the
additional structure of a brane-labeling leads to the
new combinatorial topology of the next sections.)

A window $w\in W$ on a windowed surface $F$ brane-labeled by
$\beta$ is called {\it closed} if the endpoints of $w$ coincide at
the point $p\in\delta$ and $\beta (p)=\emptyset$; otherwise, the
window $w$ is called {\it open}.

To finally explain the string phenomenology, consider a weighted
arc family in a windowed surface $F$ with brane-labeling $\beta$.
To each arc $a$ in the arc family, associate a rectangle $R_a$ of
width given by the weight on $a$, where $R_a$ is foliated by
horizontal lines as before. We shall typically dissolve the
distinction between a weighted arc $a$ and the foliated rectangle
$R_a$, thinking of $R_a$ as a ``band'' of arcs parallel to $a$
whose width is the weight. Disjointly embed each $R_a$ in $F$ with
its vertical sides in $W$ so that each leaf of its foliation is
homotopic to $a$ rel $\delta$.  Taken together, these rectangles
produce a measured foliation of a closed subsurface of $F$ as
before, and the leaves of the corresponding unmeasured vertical
foliation represent the strings.

Thus, a weighted arc family in a brane-labeled windowed surface
represents a string interaction.  Given such surfaces $F_i$ with
weighted arc families $\alpha _i$ and a choice of window $w_i$ of
$F_i$, for $i=1,2$, suppose that the sum of the weights of the
arcs in $\alpha _1$ meeting $w_1$ agrees with the sum of the
weights of the arcs in $\alpha _2$ meeting $w_2$.  In this case as
in open/closed cobordism (see e.g.\ \cite{LP}), we may glue the
surfaces $F_1,F_2$ along their windows $w_1,w_2$ respecting the
orientations so as to produce another oriented surface $F_3$, and
because of the condition on the weights, we can furthermore
combine $\alpha _1$ and $\alpha _2$ to produce a weighted arc
family $\alpha _3$ in $F_3$ (cf. Figure~5).  This describes the
basic gluing operations, namely, the operations of a c/o structure
on the space of all weighted arc families in brane-labeled
windowed surfaces (cf. Section~2 for full details). Furthermore,
these operations descend to the chain and homology levels as well
(cf. Section~2 and Appendix A).

As we shall explain (in Section~3), the degree zero indecomposables of the c/o structure are illustrated in Figures~3
and 4, and further useful degree one indecomposables are illustrated in Figure~6 (whose respective parts a-e
correspond
to those of Figure~1.)

Relations in the c/o structure of
weighted arc families or measured foliations are derived from decomposable elements, i.e., from the
fact that a given surface admits many different decompositions into ``generalized pairs of pants''
(see the next section), so
the weighted arc families or measured foliations in it can be described by different
compositions of indecomposables in the c/o structure.

We shall see that all of the known equations of open/closed string
theory, including the ``commutative and symmetric Frobenius algebras, Gerstenhaber-Batalin-Vilkovisky, Cardy,
and center (or knowledge)'' equations, hold for
the c/o structure on chains on weighted arc families (cf. Figures
7-11).

Furthermore, we shall derive several new such equations (cf. Figure 12) and in particular a set of
four new equations which together with known relations generate closed/open string duality (see
Theorem \ref{thmtrans}).

Indeed, it is relatively easy to generate many new equations of string interactions in this way, and we shall
furthermore (in Section~3.2) describe an algorithm for generating all equations of all degrees on the topological
level,
and in a sense also on the chain level.

We turn in Section 4 to the algebraic analysis of Section 3 and
derive independent sets of generators and relations in degree zero
on the topological, chain, or homology levels. In particular, this
gives a new non-Morse theoretic  calculation of the open/closed
cobordism group in dimension two \cite{Lew,LP}.  Several results
on higher degree generators and relations are also presented, and
there is furthermore a description of algebras over our c/o
structure on arc families.

Having completed this ``tour'' of the figures and this general physical discussion of the discoveries and
results contained in this paper, let us next state an ``omnibus'' theorem likewise intended to summarize the
results mathematically:

\vskip .2in

\noindent {\bf Theorem}~\it
For every brane-labeled windowed surface $(F,\beta )$, there is a space $\A(F,\beta )$ of mapping class group orbits
of suitable measured foliations in $F$ together with geometrically natural operations of gluing surfaces and measured
foliations along windows.  These operations descend to the level of piecewise-linear or cubical chains for example.

These operations furthermore descend to the level of integral homology and induce the structure of a modular
bi-operad, cf. \cite {LP}.  Algebras over this bi-operad satisfy the expected equations
as articulated in Theorem~\ref{homthm}.

Furthermore, new equations can also be derived in the language of
combinatorial topology: pairs of ``generalized pants decompositions'' of a common brane-labeled
windowed surface give rise to families of relations.

In degree zero on the homology level, we rederive the known presentation of
the open/closed cobordism groups \cite{Lew,LP}, and further partial algebraic results
are given in higher degrees.
In particular, several new relations
(which have known transformation laws) are shown to act transitively on the set of all generalized pants
decompositions of a fixed brane-labeled windowed surface.  \rm

\vskip .2in

This paper is organized as follows.  $\S$1 covers the basic combinatorial topology of measured foliations in
brane-labeled windowed surfaces and their generalized pants decompositions leading up to a description of the
indecomposables of our theory, which in a sense go back to the 1930's.  $\S$2 continues in a similar
spirit to combinatorially define the spaces $\A(n,m)$ underlying our algebraic structure on the topological
level as well as the
basic gluing operations on the topological level.  The
operations on the chain level then follow tautologically.  The operations on the homology level require the analysis
of certain fairly elaborate flows, which are defined and studied in Appendix~C and also discussed in $\S$2.  In $\S$3
continuing with combinatorial topology, we present generators, relations, and finally prove the result
that appropriate moves act transitively on generalized pants decompositions.  $\S$4 finally turns to the algebraic
discussion of the material described in $\S$3 and explains the precise sense in which the figures actually represent
traditional algebraic equations; $\S$4 furthermore presents our new algebraic results about string theory.
Closing remarks in particular include a discussion of how one might imagine our results extending
from topological to conformal field theory.

Appendix A gives the formal algebraic definition and basic
properties of a c/o structure, and  Appendix
B briefly surveys Thurston's theory of measured foliations from
the 1970-80's and describes the extension of the current paper to
the setting of general measured foliations on windowed surfaces.
It is fair to say that Appendix~A could be more appealing to a mathematician than a physicist (for whom we have tried
to make Appendix~A optional by emphasizing the combinatorial topology in the body of the paper), and that the
physically speculative Appendix~B should probably be omitted on a first reading in any case.

Appendix~C defines and studies certain flows which are fundamental
to the descent to homology as described in Appendix~A.
Nevertheless, the discussion of the flows and their salient
properties in Appendix~C is independent of the technical aspects
of Appendix~A (since chains are interpreted simply as
parameterized families); in a real sense, Appendix~C is the
substance of this paper beyond the combinatorial topology,
algebraic structure, and phenomenology, so we have strived to keep
it generally accessible.

\vskip .3in

\section{Weighted arc families, brane-labeling, and generalized pants decompositions}

\vskip .2in

\subsection{Weighted arc families in brane-labeled windowed surfaces.}

\vskip .2in

In the notation of the introduction, consider a windowed surface
$F=F_g^s(\delta _1,\ldots ,\delta _r)$, with punctures $\sigma$, boundary
distinguished points $\delta=\delta _1\cup\cdots\cup\delta _r$ and
windows $W$, together with a brane-labeling $\beta
:\delta\cup\sigma\to{\mathcal P}(B)$. Define the sets $$\delta
(\beta )=\{ p\in\delta :\beta (p)\neq\emptyset\},$$ $$ \sigma
(\beta )=\{ p\in\sigma :\beta (p)\neq\emptyset\} .$$

Fix some
brane label $A\in{\mathcal B}$, and define the brane-labeling $\beta _A$ to  be the
constant function on $\delta\cup\sigma$ with value $\{ A\}$; $\beta_A$ corresponds to the ``purely open
sector with a space-filling brane-label''.  On the other hand, the constant function $\beta _\emptyset$
with value
$\emptyset$ corresponds to the ``purely closed sector''.

It is also useful to have the notation $F^s_{g,(\# \delta_1,\ldots ,\# \delta _r)}$, where $\# S$ is the cardinality
of a set $S$.
For instance, a pair of pants with one distinguished point on each boundary component is a surface of type
$F^0_{0,(1,1,1)}$, while
the data of the windowed surface
$F_g^s(\delta _1,\delta _2,\delta _3)$ includes the specification of one point in each boundary component as well.
One further point
of convenient notation is that we shall let simply $F_{g,r}^s$ denote a surface of genus $g$ with $s$ punctures
and $r>0$ boundary
components when there is a unique distinguished point on each boundary component.

Define a $\beta$-arc $a$ in $F$ to be an arc properly embedded in $F$ with its endpoints in $W$ so that $a$ is not
homotopic fixing its endpoints into $\partial F-\delta (\beta )$.  For example, given a distinguished point
$p\in\partial
F$, consider the arc lying in a small neighborhood that simply connects one side of $p$ to another in $F$; $a$ is
a $\beta$-arc if and only if
$\beta (p)\neq\emptyset$.

Two $\beta$-arcs are {\it parallel} if they are homotopic rel $\delta$, and a $\beta$-arc family is the homotopy
class
rel $\delta$ of a collection of $\beta$-arcs, no two of which are parallel.  Notice that we take homotopies rel
$\delta$
rather than
rel $\delta (\beta )$.

A {\it  weighting} on an arc family is the assignment of
a positive  real number to each of its components.

Let $Arc'(F,\beta )$ denote the geometric realization of the partially
ordered set of all $\beta$-arc families in $F$.  $Arc'(F,\beta )$ is
described as the set of all projective positively weighted $\beta$-arc
families in $F$ with the natural topology. (See for instance \cite{KLP} or \cite{P1} for
further details and Figure~2 below for an illustrative example.)

The {\it (pure) mapping class group} $MC(F)$ of $F$ is the group of
orientation-preserving homeomorphisms of $F$ pointwise fixing
$\delta\cup\sigma$ modulo homotopies pointwise fixing $\delta\cup\sigma$.  $MC(F)$
acts naturally on $Arc'(F,\beta )$ by definition with quotient the {\it arc
complex}
$$Arc(F,\beta )=Arc'(F,\beta )/MC(F).$$

We shall also require the corresponding deprojectivized versions:
$\widetilde{Arc}'(F,\beta )\approx Arc'(F,\beta )\times{\mathbb R}_{>0}$ is the space of all positively weighted arc
families in $F$ with the natural topology, and $$\widetilde{Arc}(F,\beta )=\widetilde{Arc}'(F,\beta )/MC(F)\approx
Arc(F,\beta )\times{\mathbb R}_{>0}.$$

\vskip .2in

\hskip 1.3in\epsffile{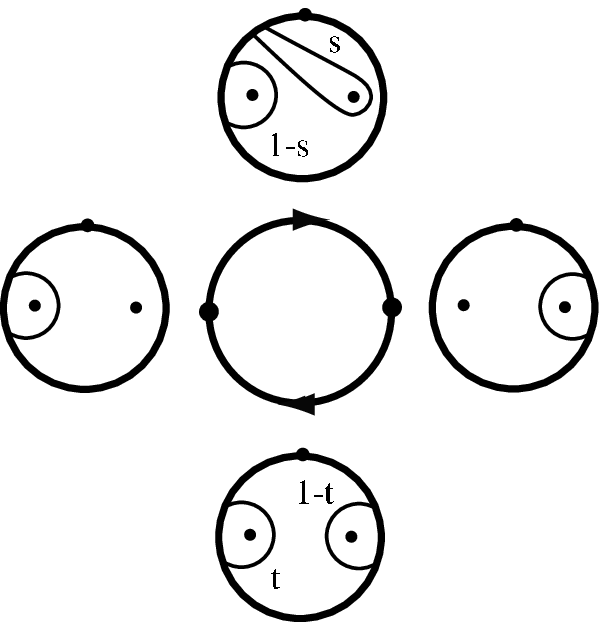}

\vskip .2in

\noindent \begin{footnotesize}
{\bf Figure 2}~The arc complex $Arc(F_{0,1}^2,\beta _\emptyset )$ is homeomorphic to a circle $S^1$.
We omit the common label
$\emptyset$ at each point of $\delta\cup\sigma$ to avoid cluttering the figure.  There are exactly the two
$MC(F_{0,1}^2)$-orbits of $\beta _\emptyset$-arcs on the right and left.  These can be disjointly embedded
in the two distinct ways at the top and bottom.  As the parameter $t$ on the bottom
varies in the range $0\leq t\leq 1$, there is described a projectively weighted $\beta _\emptyset$-arc family,
that is, the
two disjoint arcs determine a one-dimensional simplex in $Arc(F,\beta _\emptyset )$, and likewise for the
parameter $s$ on
the top. The two one-simplices are incident at their endpoints as illustrated to form a circle ${Arc}(F,\beta
_\emptyset)\approx S^1$.  Furthermore, $Arc'(F,\beta _\emptyset )\approx{\mathbb R}$, $\widetilde{Arc}(F,\beta
_\emptyset )\approx S^1\times{\mathbb R}_{>0}$, and $\widetilde{Arc}'(F,\beta _\emptyset )\approx
{\mathbb R}\times{\mathbb
R}_{>0}$ with the primitive mapping classes acting by translation by one on ${\mathbb R}$.
\end{footnotesize}

\vskip .2in

It will be useful in the sequel to employ a notation similar to
that in Figure~2, where parameterized collections of arc families
are described by pictures of arc families together with functions
next to the components, where the functions represent the
parameterized evolution of weights.  We shall also typically let
the icon $\bullet$ denote either a puncture or a distinguished
point on the boundary as in Figure~2.

In contrast to Figure~2, if we instead consider the purely open
sector with space-filling brane-label $\beta _A$, for
$A\neq\emptyset$, then there is yet another $MC(F_{0,1}^2)$-orbit
of arc encircling the boundary distinguished point.  In this case,
$Arc(F_{0,1}^2,\beta _A)$ is homeomorphic to the join of the
circle in Figure~2 with the point representing this arc, namely,
$Arc(F_{0,1}^2,\beta _A)$ is homeomorphic to a two-dimensional
disk.

For another example of an arc complex, take the brane-labeling
$\beta _\emptyset\equiv\emptyset$ on $F_{0,2}^0$,
for which again $MC(F_{0,2}^0)\approx{\mathbb Z}$.
There is a unique $MC(F_{0,2}^0)$-orbit of singleton $\beta _\emptyset$-arc, and there are two possible
$MC(F_{0,2}^0)$-orbits of $\beta
_\emptyset$-arc families with two component arcs illustrated in Figure~3.  Again, $Arc(F_{0,2}^0,\beta _\emptyset )$
is
homeomorphic to a circle.  If $A\neq\emptyset$, then $Arc(F_{0,2}^0,\beta _A)$ is homeomorphic to a three-dimensional
disk.

To explain the connection  with earlier
work, consider the purely closed sector
$\beta_\emptyset\equiv\emptyset$ on $F=F_{g,r}^s$. Let
$Arc_{\#}(F)$ denote the subspace of $Arc(F,\beta _\emptyset)$
corresponding to all projective positively weighted arc families
$\alpha$ so that each component of $F-\cup\alpha$ is either a
polygon or an exactly once-punctured polygon, i.e., $\alpha$ quasi fills $F$;  $Arc_{\#}(F)$ was
shown in \cite{P2} to be proper homotopy equivalent to a natural
bundle over Riemann's moduli space of the bordered surface $F$ as defined in the Introduction provided $F$ is not an
annulus
($g=s=0$,
$r=2$). 

Let $Arc(F)\supseteq Arc_{\#}(F)$ denote the subspace of
$Arc(F,\beta _\emptyset )$ corresponding to all projective
positively weighted arc families $\alpha$ so that each window of
$F$ (i.e., each boundary component) has at least one arc in
$\alpha$ incident upon it. The spaces $Arc(F)$ comprise the
objects of the basic topological operad studied in \cite{KLP}.

\vskip .2in

\subsection{Generalized pants decompositions.}
\label{indec}

A {\it generalized pair of pants} is a surface of genus zero with $r$ boundary components and $s$
punctures, where $r+s=3$, with exactly one distinguished point on each boundary component, that is,
a surface of type $F^0_{0,3}$, $F^1_{0,
2}$, or $F^2_{0,1}$.

A {\it (standard) pants decomposition} $\Pi$ of a windowed surface $F=F_g^s(\delta _1,\ldots ,\delta _r)$ is
(the homotopy class of) a collection of disjointly embedded essential curves in the interior $F$, no two of which are
homotopic, together with a condition on the complementary regions to $\Pi$ in $F$.

To articulate this condition, let us enumerate the curves $c_1,\ldots ,c_K$ in $\Pi$, choose disjoint annular
neighborhoods $U_k$ of $c_k$ in $F$, for $k=1,\ldots ,K$, and set $U=U_1\cup\cdots\cup U_K$.
Just for the purposes of articulation, let us also choose on each boundary component of $U$ a distinguished point.
We require that each component of $F-U$ is a generalized pair of pants
or a boundary-parallel annulus of type $F^0_{0,(1,n)}$, for some $n\geq 1$.

Simple Euler characteristic considerations give the
following lemma.

\vskip .2in

\begin{lem}\label{eulerlemma} For a windowed surface $F=F_g^s(\delta _1,\ldots ,\delta _r)$, there are
$\#W=\sum _i^r \#\delta _i=\#\delta$ many windows.
The real dimension of $Arc(F,\beta _\emptyset )$ is $6g-7+3r+2s+\#\delta$.
Furthermore, there are $3g-3+2r+s$ curves in a pants
decomposition $\Pi$ of $F$ and $2g-2+r+s$ generalized pairs of pants complementary to an annular
neighborhood of the pants curves.\end{lem}

\vskip .2in

If $\beta$ is a brane-labeling on the windowed surface $F$, then a {\it generalized pants decomposition} of
$(F,\beta )$ is (the homotopy class of) a family of disjointly embedded closed curves in the interior of $F$ and arcs
with endpoints
in $\delta (\beta )\cup\sigma (\beta)$, no two of which are parallel,  so that each complementary region is one
of the
following {\it indecomposable brane-labeled surfaces}:

\begin{itemize}

\item[--]
a triangle $F^0_{0,(3)}$ with no vertex brane-labeled by $\emptyset$;

\item[--] a generalized pair of pants $F^0_{0,3}$, $F^1_{0,
2}$, or $F^2_{0,1}$ with all points $\delta$ in the boundary brane-labeled by $\emptyset$;

\item[--]
 a once-punctured monogon $F^1_{0,1}$ with
puncture brane-labeled by $\emptyset$ and boundary distinguished point by $A\neq\emptyset$;

\item[--]

an annulus  $F^0_{0,2}$
with at least point of $\delta$ labeled by $\emptyset$.
\end{itemize}

For instance, if every brane label is empty, then a generalized pants decomposition is a standard pants
decomposition.   At the other extreme, if every brane label is non-empty, then $F$ admits a decomposition into
triangles
and
once-punctured monogons, a so-called {\it quasi triangulation}  of $F$, cf. \cite{P3}; see Figures~12b and 13 for
examples. Provided there is at least one non-empty brane
label, we may collapse each boundary component with empty brane label to a puncture to produce another windowed
surface
$F'$ from
$F$.  A quasi triangulation of $F'$ can be completed with brane-labeled annuli to finally produce a generalized
pants
decomposition
of $F$ itself.

Thus, any brane-labeled windowed surface admits a generalized pants decomposition.  Furthermore, any collection of
disjointly embedded essential curves and arcs connecting non-empty brane labels so that no two components are
parallel can
be completed to  a generalized
pants decomposition.

\vskip .2in

\subsection{Indecomposables.}

\vskip .2in

We shall introduce standard foliations on indecomposable surfaces which are the basic
building blocks of the theory, and we begin with the annulus in Figure~3.

\vskip .2in

\hskip .4in\epsffile{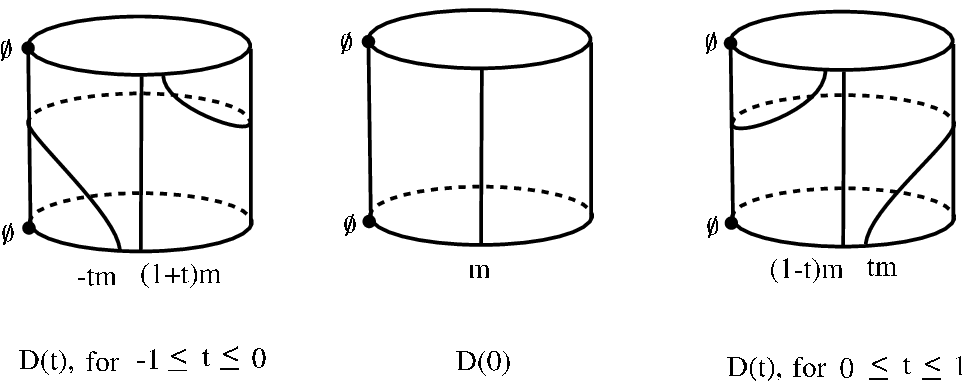}

\centerline{ \begin{footnotesize}{\bf Figure 3} The twist flow on
$\widetilde{Arc}'({\mathcal A},\beta _\emptyset)$.\end{footnotesize}}

\vskip .2in

In the notation of Figure~3, consider the purely closed sector with
 brane-labeling $\beta_\emptyset\equiv\emptyset$ on
a fixed annulus ${\mathcal A}$ of type $F_{0,2}^0$. Define a
one-parameter ``Dehn twist flow'' $D(t)$, for $-1\leq t\leq 1$, on
$\widetilde{Arc}'({\mathcal A},\beta _\emptyset )$, as illustrated
in the figure, where $m$ denotes the sum of the weights of the
arcs in $\alpha\in \widetilde{Arc}'({\mathcal A},\beta _\emptyset
)$.  Letting $T$ denote the right Dehn twist along the core of the
annulus, one extends to all positive real values of $t$ by setting
$D(t)(\alpha )=T ^{[t]} D(t-[t])(\alpha )$, where $[t]$ denotes
the integral part of $t$, and likewise for negative real values of
$t$.

Figure~4 illustrates the remaining building blocks of the theory.
Notice that $F_{0,1}^1$ brane-labeled with some $\beta$ taking value $\emptyset$ on the boundary is absent from
Figure~4 and implicitly from the theory since $Arc(F_{0,1}^1,\beta )$ is empty.

A fact going back to Max Dehn in the 1930's is that ``free''
homotopy classes rel $\delta (\beta _\emptyset )=\emptyset$ in a
fixed pair of pants ${\mathcal P}$ of type $F_{0,3}^0$ are
determined by the three ``intersection numbers'' $m_1,m_2,m_3$,
namely, the number of endpoints of component arcs in each
respective boundary component, subject to the unique constraint
that $m_1+m_2+m_3$ is even.  Two representative cases are
illustrated in Figure~4e, and the full partially ordered set is
illustrated in Figure~4d. There are conventions in the pair of
pants that have been suppressed here insofar as the ``arc
connecting a boundary component to itself goes around the right
leg of the pants''; see Figure~4f and see \cite{Pprobe} for
details.

\vskip .2in

\hskip .15in\epsffile{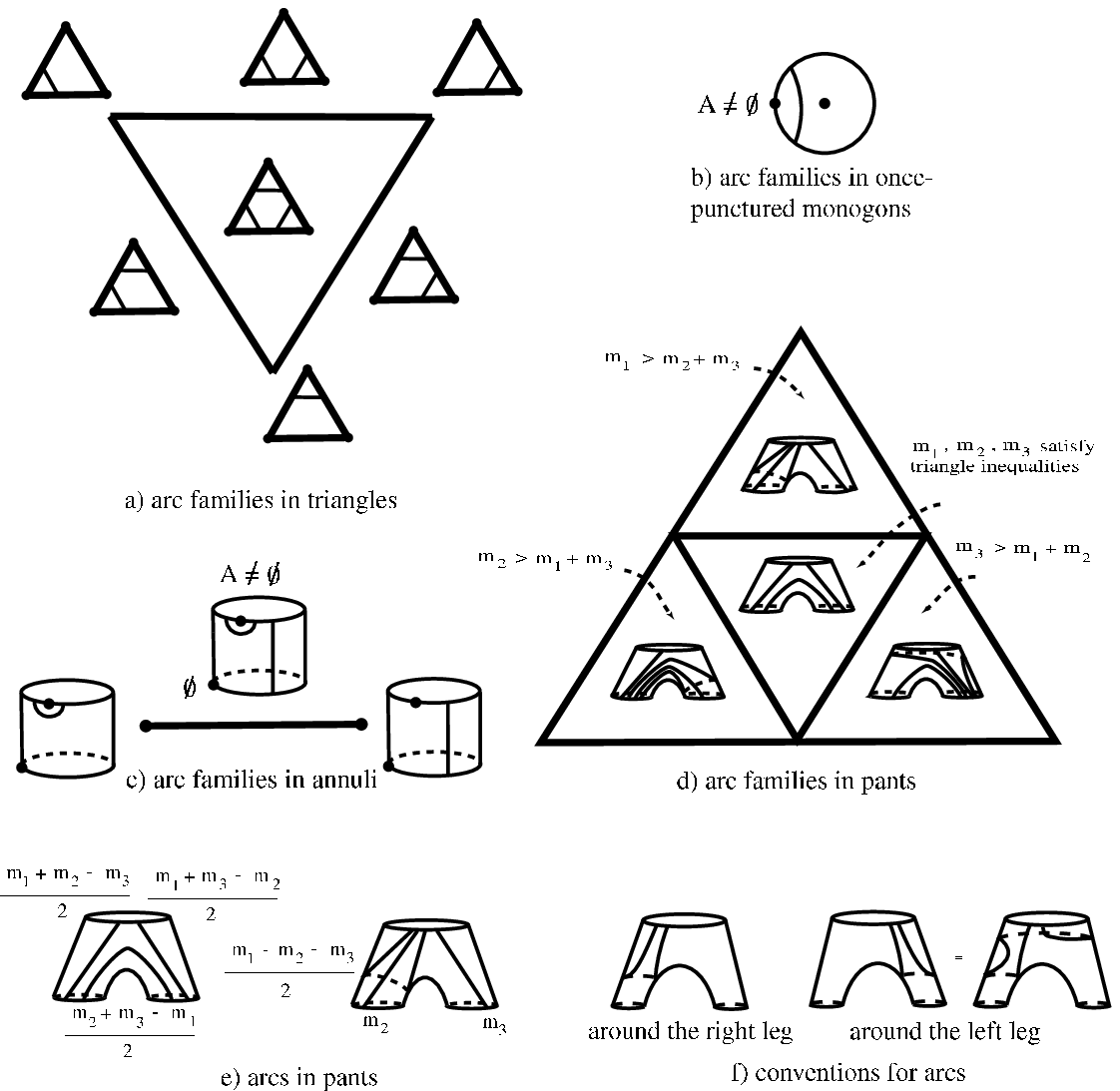}

\noindent\begin{footnotesize}{{\bf Figure 4} Indecomposables.}  We depict the geometric realization of
$Arc(F^0_{0,(3)},\beta )$ in part a for some brane-labeling
$\beta$ whose image does not contain $\emptyset$, and which is simply omitted from the figure.
There is the unique
element of $Arc(F_{0,1}^1,\beta )$ depicted in part b when the brane-label on the boundary is non-empty.
For the brane-labeling $\beta$ on $F_{0,2}^0$ indicated in part c, we consider instead the
homotopy classes of $\beta$-arc families rel $\delta (\beta )$, rather than rel
$\delta$ as before.  Likewise, for the brane-labeling
$\beta _\emptyset$ on $F_{0,3}^0$, which we omit from the figure, we consider again the homotopy classes
of $\beta_\emptyset$-arc
families rel $\delta (\beta _\emptyset)$, where $\delta (\beta _\emptyset )=\emptyset$ by
definition, and depict the geometric realization in part d.
\end{footnotesize}

\vskip .2in

One further remark is that arc families in all generalized pairs of pants are also implicitly described
by Figures~4d-f,
where punctures correspond to boundary components with no incident arcs.

\vskip .2in

\subsection{Standard models of arc families.}

\vskip .2in

Suppose that $\Pi$ is a generalized pants decomposition of a brane-labeled windowed surface $(F,\beta )$,
where $\Pi$ has curve components $c_1,\ldots c_K$ and arc components $d_1,\ldots d_L$.
Let $U_k$ denote a fixed annular
neighborhood of
$c_k$ for $k=1,\ldots ,K$, and set $U=U_1\cup\cdots\cup U_K$.

In order to parameterize weighted arc families, we must make
several further choices, as follows. Choose a framing to the
normal bundle to each curve $c_k$, which thus determines an
identification of the unit normal bundle to $c_k$ in $F$ with the
standard annulus ${\mathcal A}$.  In turn, the unit normal bundle
is also identified with the neighborhood $U_k$, and there is thus
an identification of $U_k$ with ${\mathcal A}$ determined by the
framing on $c_k$. Furthermore, choose homeomorphisms of each
generalized pair of pants component of $F-U$ with some standard
generalized pair of pants ${\mathcal P}$.  Choose an embedded
essential arc $a_0$ once and for all in ${\cal A}$, and likewise
choose standard models for arc families in ${\mathcal P}$ (say,
with the conventions for twisting as in Figure~4f). Let us call a
generalized pants decomposition together with this specification
of further data a {\it basis} for arc families.

Given $\alpha\in\widetilde{Arc}'(F,\beta )$, choose a representative weighted arc
family that meets each component of $\Pi$ transversely a minimal number of times, let $m_k$ denote the sum of the
weights
of the arcs in $\alpha$ that meet $c_k$ counted with multiplicity (and without a sign), and let $n_\ell$ denote the
analogous sum for the arcs $d_\ell$.

\vskip .2in

\begin{thm}\label{DT}
Fix a basis for arc families with underlying generalized pants decomposition $\Pi$ of a brane-labeled windowed
surface
$(F,\beta)$, and adopt the notation above given some $\alpha\in\widetilde{Arc}'(F,\beta )$.
 Then under the identifications with the standard annulus ${\mathcal A}$ and standard pants
${\mathcal P}$,
$\alpha$ is represented by a weighted arc family that meets complementary regions to $U$ in $F$ in
exactly one of the configurations shown in Figure~4 and meets each $U_k$ in $D(t_k)(a_0)$ for some well-defined
$t_k\in{\mathbb R}$, where $a_0$ is weighted by $m_k$.
Furthermore, a point of $\widetilde{Arc}'(F)$ is uniquely determined by its coordinates $(m_k,t_k)$ for $k=1\ldots K$
and $n_\ell$, for $\ell =1,\ldots ,L$.
\end{thm}

\vskip .2in

\begin{proof}
Since the arcs in a generalized pants decomposition connect points of $\delta\cup\sigma$ and the components of an arc
family avoid a
neighborhood of $\delta\cup\sigma$, intersections with triangles and once-punctured monogons are established.
We may homotope an arc family to a standard model in each pair of pants; the twisting numbers $t_k$ are then the
weighted algebraic intersection numbers (with a sign) with $a_0$ in each annulus (all arcs oriented from
top-to-bottom
or bottom-to-top of the annulus);
see   the ``Dehn-Thurston'' coordinates from
\cite{Pthesis},
\cite{Pprobe} for further details.
\end{proof}

\vskip .2in

\begin{cor} \label{DTcor}In the notation of Theorem \ref{DT}, any parameterized family in $\widetilde{Arc}'(F,\beta )$
is represented by one that meets complementary regions to $U$ in
$F$ in parameterized families of the configurations shown in
Figure~4 and meets each $U_k$ in $D(t_k)(a_0)$, where $t_k$
depends upon the parameters, for $k=1,\ldots ,K$.  Furthermore, a
parameterized family is uniquely determined by its parameterized
coordinates $(m_k,t_k)$ for $k=1\ldots K$ and $n_\ell$, for $\ell
=1,\ldots ,L$.\hfill$\Box$
\end{cor}

\vskip .2in

Notice that in either case of the theorem or the corollary, the intersection numbers on any triangle satisfy all
three possible weak triangle inequalities.

\vskip .3in

\section{C/O string operations on weighted arc families}\label{oparcs}

\vskip .2in

Recall that a window $w\in W$ in a brane-labeled windowed surface $F$ is {\it closed} if its closure is an
entire boundary
component of $F$
and the distinguished point complementary to $w$ is brane-labeled by $\emptyset$, and otherwise the window is
{\it open}.

Given a positively weighted arc family in $F$, let us furthermore say that a window $w\in W$ is {\it active}
if there
is an arc in the family with an endpoint in $w$, and otherwise the window is {\it inactive}.

In order to most directly connect with the usual phenomenology of strings, we shall require all windows to be active,
but
the more  general case of operations on inactive windows is not uninteresting, specializes to the treatment here,
and
will be discussed in Appendix~B.

Given a positively weighted arc family in $F$, we may simply
collapse each inactive window, or consecutive sequence of inactive
windows in a boundary component, to a new distinguished point on
the boundary, where the brane-labeling of the resulting
distinguished point is the union of all the brane labels on the
endpoints of the windows collapsed to it.  In case a boundary component
consists entirely of inactive windows, then it is collapsed to a
new puncture, which is again brane-labeled by the union of all the
brane labels on the collapsed boundary component.
Thus, given any positively weighted arc family in $F$,
there is a corresponding positively weighted arc family in a
corresponding surface so that each window is active.
(This is one explanation
for why we brane-label by the power-set of branes, namely, in order to
effectively
take every window to be active.)

For any windowed surface $F$, define
$$\widetilde{Arc}(F)=\bigsqcup \widetilde{Arc}(F,\beta),$$
where the disjoint union is over all brane-labelings on $F$.
The basic objects of our topological c/o structure are
$$\widetilde{Arc}(n,m)=\bigsqcup \left\{ \alpha\in \widetilde{Arc}(F):
\begin{array}{l}\alpha ~{\rm has}~n~{\rm closed~and}~ m~{\rm
open~active}\\{\rm
windows~and~no~inactive~windows}\end{array}\right\},$$ where the
disjoint union is over all orientation-preserving
homeomorphism classes of windowed surfaces.

\vskip .2in

\centerline{\epsffile{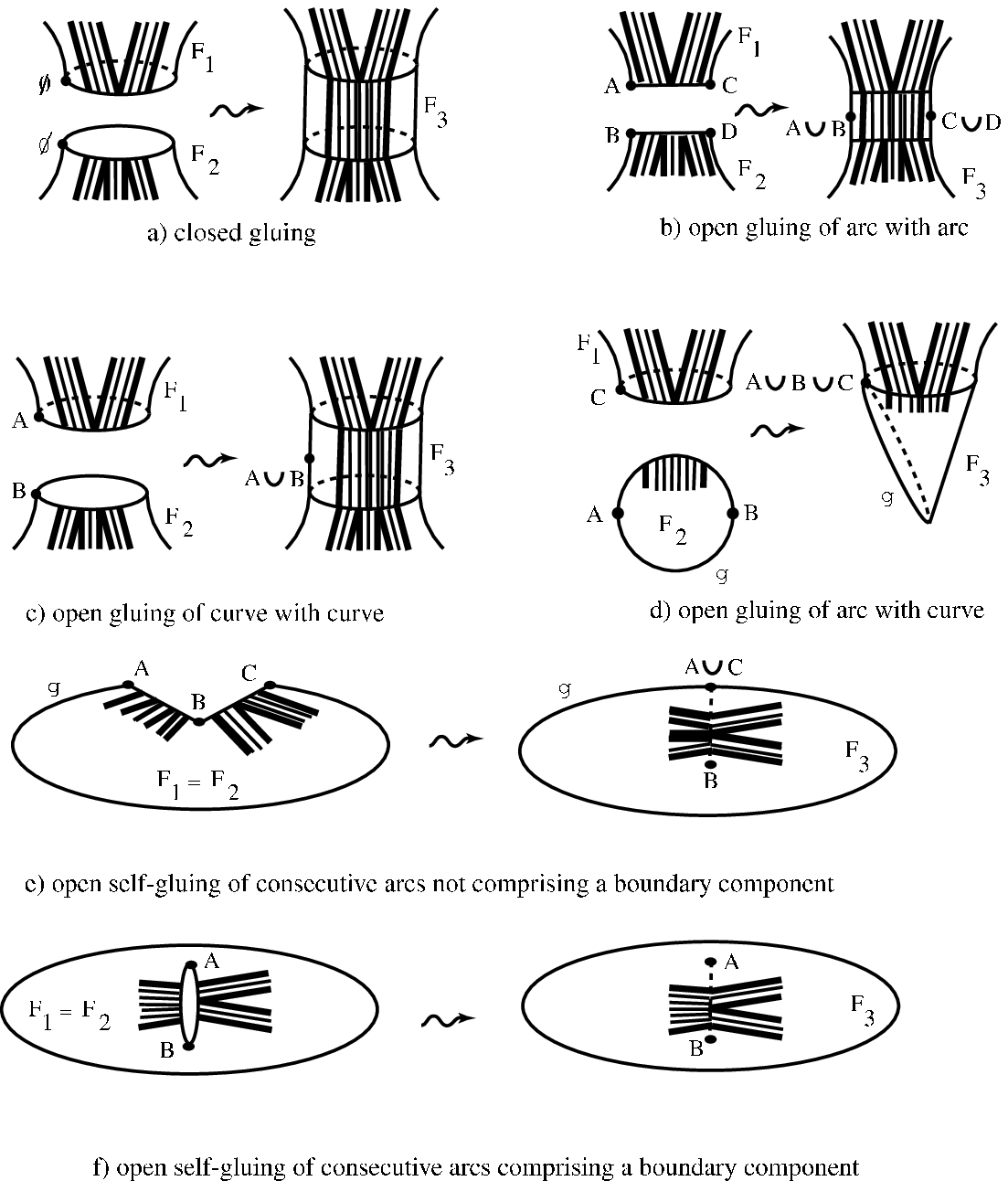}}

\vskip .2in
\centerline{\begin{footnotesize}{\bf Figure 5} The c/o operations on measured foliations.\end{footnotesize}}

\vskip .2in

If $\alpha\in\widetilde{Arc} (n,m)$, then define the {\it $\alpha$-weighting} of an active window $w$ to be the sum
of
the weights of arcs in $\alpha$ with endpoints in $w$, where we count with multiplicity
(so if an arc in $\alpha$ has both endpoints in $w$, then the weight of this arc contributes twice to the weight of
$w$).

Suppose we have a pair of arc families $\alpha_1,\alpha _2$
in respective windowed surfaces $F_1,F_2$ and a pair of active windows $w_1$ in $F_1$ and $w_2$ in $F_2$, so that the
$\alpha _1$-weight of $w_1$ agrees with the $\alpha _2$-weight of $w_2$.
Since $F_1,F_2$ are oriented surfaces, so too are the windows $w_1,w_2$ oriented.   In each operation, we
identify windows reversing orientation, and we identify certain distinguished points.

To define the open and closed gluing ($F_1\neq F_2$) and self-gluing ($F_1=F_2$) of $\alpha_1,\alpha _2$ along the
windows
$w_1,w_2$, we
identify windows and distinguished points in the natural way and combine foliations.  In closed string operations,
we
``replace the  distinguished point, so there is no puncture'' whereas with open string operations, ``distinguished
points
always beget either other distinguished points or perhaps punctures''.  In any case whenever distinguished points
are
identified, one takes the union of brane labels (the intersection of branes) at the new resulting distinguished
point or
puncture.

More explicitly, the general procedure of gluing defined above
specializes to  the following specific operations on the $\widetilde{Arc}
(n,m)$:

\vskip .1in

\noindent {\bf Closed gluing and self-gluing}~See Figure 5a.
Identify the two corresponding boundary components of $F_1$ and
$F_2$, identifying also the distinguished points on them {\sl and
then including this point in the resulting surface} $F_3$. $F_3$
inherits a brane-labeling from those on $F_1,F_2$ in the natural
way.  We furthermore glue $\alpha _1$ and $\alpha _2$ together in
the natural way, where the two collections of foliated rectangles
in $F_1$ and $F_2$ which meet $w_1$ and $w_2$ have the same total
width by hypothesis and therefore glue together naturally to
provide a measured foliation ${\mathcal F}$ of a closed subsurface
of $F_3$.  (The projectivization of this gluing operation is
precisely the composition in the cyclic operad studied in
\cite{KLP}; we have deprojectivized and included the weighting
condition in the current paper in order to allow self-gluing of
closed strings as well.)

\vskip .1in

\noindent {\bf Open gluing}~The surfaces $F_1$ and $F_2$ are distinct, and we identify $w_1$ to $w_2$ to produce
$F_3$.
There are cases depending upon whether the closure of $w_1$ and $w_2$ is an interval or a circle.  The salient
cases are
illustrated in Figure 5b-d.  In each case,  distinguished points on the boundary in $F_1$ and $F_2$ are
identified to produce a new distinguished boundary point in $F_3$, and the brane labels are combined, as is also
illustrated.  As before, since the
$\alpha _1$-weight on $w_1$ agrees with the $\alpha _2$-weight on $w_2$, the foliated rectangles again combine to
provide a measured foliation ${\mathcal F}$ of a closed subsurface of $F_3$.

\vskip .1in

\noindent {\bf Open self-gluing}~There are again cases depending upon whether the closure of $w_1$ or $w_2$ is a
circle or an interval, but there is a further case as well when the two intervals lie in a common boundary component
and
are consecutive.  Other than this last case, the construction is identical to those illustrated
in Figure~5b-d.  In case the two windows are consecutive along a common boundary component, again they are identified
so
as to produce surface $F_3$ with a puncture resulting from their common endpoint as in Figure~5e-f, where the
puncture is
brane-labeled by the label of this point, and the foliated rectangles combine to provide a measured foliation
${\mathcal F}$ of a closed subsurface of $F_3$.

\vskip .1in

At this stage, we have only constructed a measured foliation ${\mathcal F}$ of a closed subsurface of $F_3$, and
indeed,
${\mathcal F}$ will typically not be a weighted arc family.  By Poincar\'e recurrence, the sub-foliation
${\mathcal
F}'$ comprised of leaves that meet $\partial F$ corresponds to a weighted arc family $\alpha _3$ in $F_3$.
Notice that
the $\alpha _3$-weight of any window uninvolved in the operation agrees with its $\alpha _1$- or $\alpha _2$-weight,
so
in particular, every window of $F_3$ is active for $\alpha _3$.

Let us already observe here that the part of  ${\mathcal F}$  that we discard to get $\alpha _3$ can naturally
be
included (as we shall discuss in Appendix~B).  Furthermore, notice that a gluing operation never produces a
``new'' puncture brane-labeled by $\emptyset$.

The assignment of $\alpha_3$ in $F_3$ to $\alpha _i$ in $F_i$, for $i=1,2$ completes the definition of the
various operations.  Associativity and equivariance for bijections are immediate, and so we have our first
non-trivial
example of a c/o structure (see Appendix~A for the precise definition):

\vskip .2in

\begin{thm}
\label{cothm}
Together with open and closed gluing and self gluing operations,
the spaces $\widetilde{Arc}(n,m)$ form a topological c/o structure.
Furthermore, this c/o structure is brane-labeled by $\mathcal{P}(\B)$ and is  a
$(g,\chi -1)$-c/o structure, where $g$ is the genus and $\chi$ is the
Euler-characteristic.
\end{thm}

\vskip .2in

\begin{proof}
See Appendix~A for the definitions and the proof.
\end{proof}

\vskip .2in

\begin{cor}
The open and closed gluing operations descend  to operations on
the PL chain complexes of $\widetilde{Arc}(n,m)$ giving them a
chain level c/o structure.
\end{cor}

\begin{proof}
We define a ``chain level c/o structure''  in such a manner that
this follows immediately from the previous theorem; see Appendix
A for details.
\end{proof}

\vskip .2in

\begin{thm}
\label{archom} The integral homology groups $H_*(\widetilde{Arc}(n,m))$ comprise a
modular bi-operad when graded by genus for closed gluings and
self-gluings and by Euler-characteristic-minus-one for open
gluings and self-gluings.
\end{thm}

\begin{proof}
In contrast to the previous corollary, this requires more than
just a convenient definition since we must first show that the
gluing operations descend to the level of homology; specifically,
given homology classes in $\widetilde{Arc}(n,m)$ and
$\widetilde{Arc}(n',m')$, we must find representative chains that
assign a common weight on the windows to be glued.

This is accomplished by
introducing two continuous flows on $\widetilde{Arc}(F,\beta )$ for
each window
$w$, namely,
$\psi _t^{w}$, for $0\leq t\leq 1$ for non self-gluing and $\phi _t^w$, for $-1\leq t<1$  for
self-gluing, where $\beta$ is a fixed brane-labeling on the windowed surface $F$.  In
effect for non self-gluing,
$\psi ^w_t$ simply scales in the ${\mathbb R}_{>0}$-action on $\widetilde{Arc}(F,\beta )$ so that the
weight of window
$w$ is unity at time one.  To describe the key attributes of the more complicated flow for self-gluing,
suppose that $w'\neq w$ is any other window of $F$ and $\alpha\in\widetilde{Arc}(F,\beta )$ where
the $\alpha$-weight of $w'$ is less than the $\alpha$-weight of $w$.

 There
is a well-defined ``critical'' value $t_c=t_c(\alpha )$ of $t$ so that the $\phi ^{w}_{t_c}(\alpha
)$-weight of $w$  first agrees with the $\phi ^{w}_{t_c}(\alpha ')$-weight of $w'$; furthermore,
the function $t_c(\alpha )$ is continuous in
$\alpha$.

These flows are defined and studied in Appendix~C, and the theorem then
follows directly from Proposition~\ref{homdescent}.
\end{proof}

\section{Operations, Relations, Duality}

\subsection{Operations}
Operations may be conveniently described by weighted arc families,
or by parameterized families of weighted arc families.  If a
parameterized family of arc families depends upon $p$ real
parameters, then we shall say that it is an operation of {\it
degree} $p$.  In order to establish notation, the standard
operations in degrees zero and one are illustrated in Figure~6.

In this figure, the distinguished points on the boundary
come with an enumeration that we have typically suppressed.  Only for clarity for the bracket in Figure~6k do we
indicate
the enumeration of the distinguished points with the numerals ``1'' and ``2''; we shall omit such enumerations
in subsequent figures since they can be inferred from the incidence and labeling of arcs in the figure.

It is worth remarking that the BV operator $\Delta _a(t)$ is none other than the projection
to $MC(F)$-orbits of the Dehn twist operator $D(t)$ discussed in Section~1.3.

\vskip .2in

\centerline{\epsffile{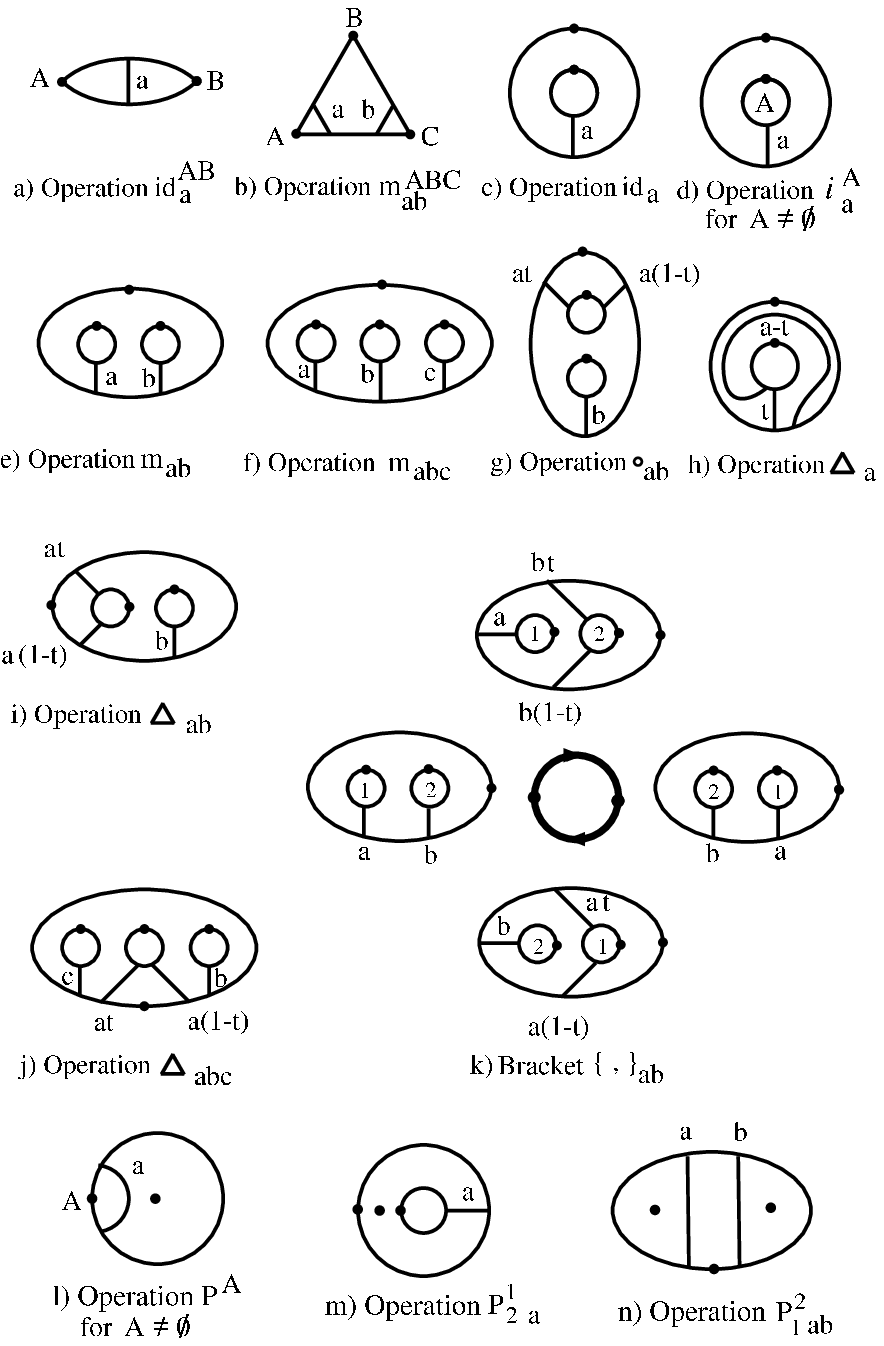} }

\begin{footnotesize}
\noindent{{\bf Figure 6}~Standard operations of degrees zero and one.
If
there is no brane-label indicated, then the label is tacitly taken to be
$\emptyset$. See Section~\ref{chainops} for the traditional algebraic interpretations.}
\end{footnotesize}

\subsection{Relations}

\vskip .2in

Relations in the c/o structure on $\widetilde{Arc}(n,m)$ or its
chain complexes can be described and derived by fixing some
decomposable windowed surface $F$, choosing two generalized pants
decompositions $\Pi ,\Pi '$ of $F$ and specifying an arc family
$\alpha$ or a parameterized family of arc families in $F$.  Each
of $\Pi $ and $\Pi '$ decompose $F$ into indecomposable surfaces
and annular neighborhoods of the pants curves.

According to Theorems~\ref{DT} and
\ref{cothm}, $\alpha$ thus admits two different descriptions as iterated compositions of
operations in the c/o structure, and these are equated to derive the corresponding algebraic relations.  We shall
abuse
notation slightly and simply write an equality of two pictures of $F$, one side of the equation illustrating $\alpha$
and $\Pi $ in $F$ and the other illustrating $\alpha$ and $\Pi '$; we shall explain the algebraic interpretations in
the
next section.  As with operations, a relation on a
$p$ parameter family of weighted arc families is said to have {\it degree} $p$.

Accordingly, Figures~7 and 8 illustrate all of the standard
relations of two-dimensional open/closed cobordism (cf.
\cite{Lew,Segal2,MS,LP}). In particular, notice that the
``Whitehead move'' in Figure~7a  corresponds to associativity of
the open string operation. The Cardy equation in Figure~7e depends
upon the two generalized pants decompositions $\Pi ,\Pi '$ of the
surface $F^0_{0,2}$  with no empty brane labels, where $\Pi$
consists of a single simple closed curve, and $\Pi '$ is an ideal
triangulation.

\vskip .1in

\centerline{\epsffile{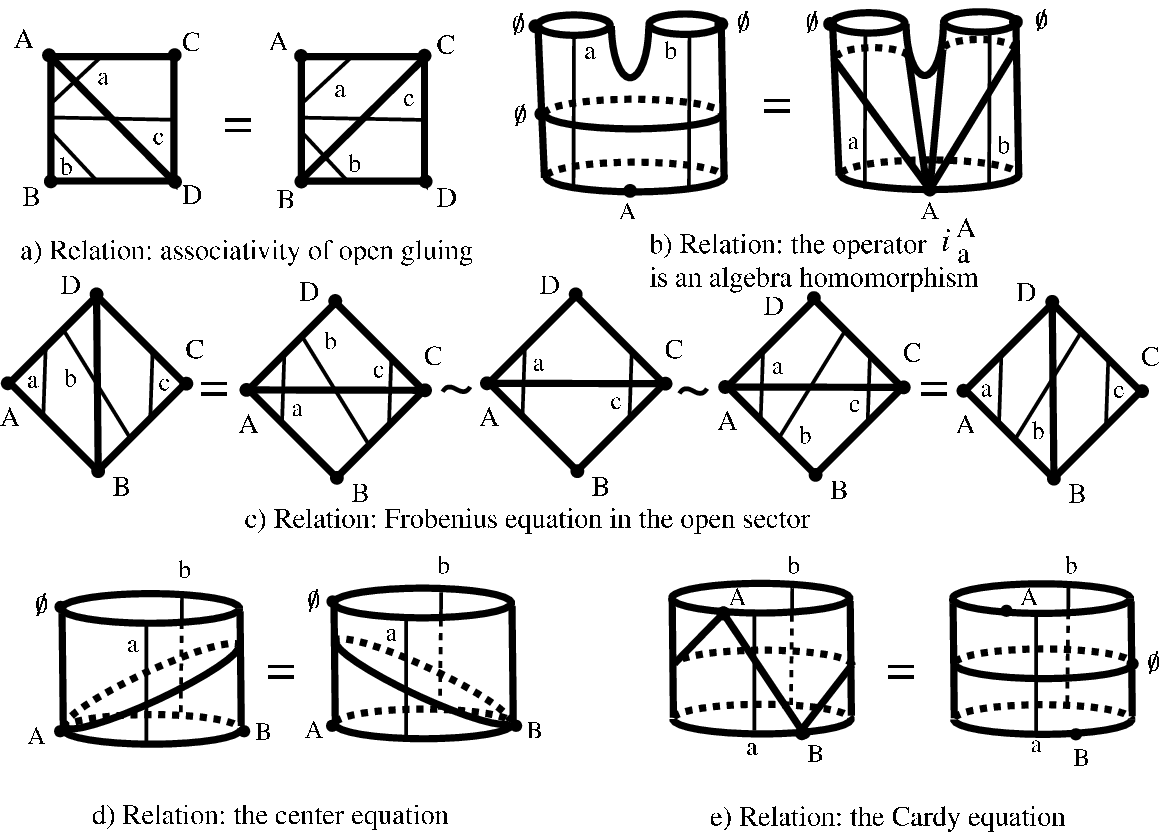}}

\begin{footnotesize}
\centerline {{\bf Figure 7} Open/closed cobordism relations.}
\end{footnotesize}

\vskip .2in

The Frobenius equation is more interesting since it consists of
two pairs $(\alpha _i, \Pi _i)$, for $i=1,2$, where $(\alpha _i
,\Pi _i)$ is comprised of a weighted arc family $\alpha_i $ with
each window active and an ideal triangulation $\Pi _i$ of a
quadrilateral; see Figure~7c at the far left and right.  Perform
the unique possible Whitehead move on $\Pi _i$ to get $\Pi '_i$,
for $i=1,2$.  In fact, the pairs $(\alpha _1,\Pi _1')$ and
$(\alpha _2,\Pi _2')$ are not identical, rather they are homotopic
in $\widetilde{Arc}(0,4)$, as is also illustrated in Figure~7c.

\vskip .2in

\centerline{\epsffile{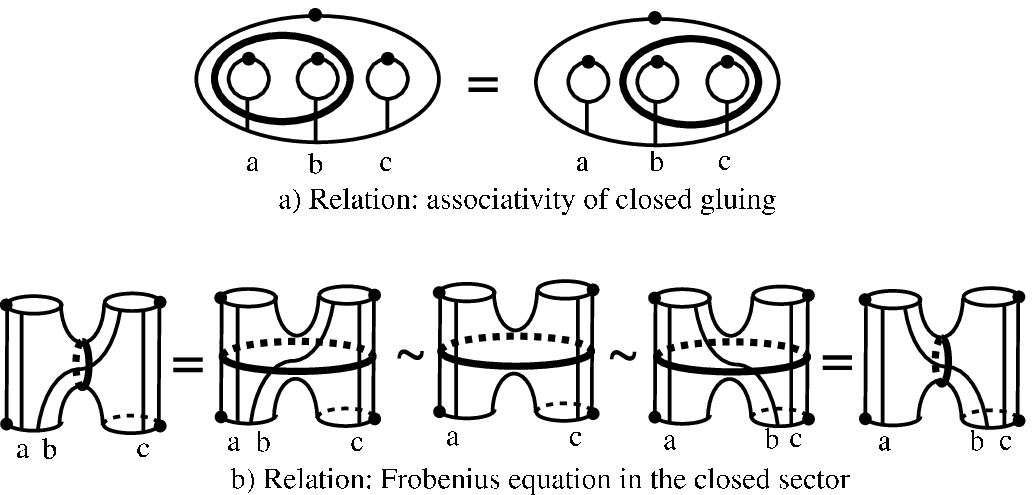}}

\vskip .2in

\begin{footnotesize}
\noindent
{{\bf Figure 8} Further open/closed cobordism relations: associativity and the Frobenius equation in the closed
sector.}
\end{footnotesize}

\vskip .2in

\centerline{\epsffile{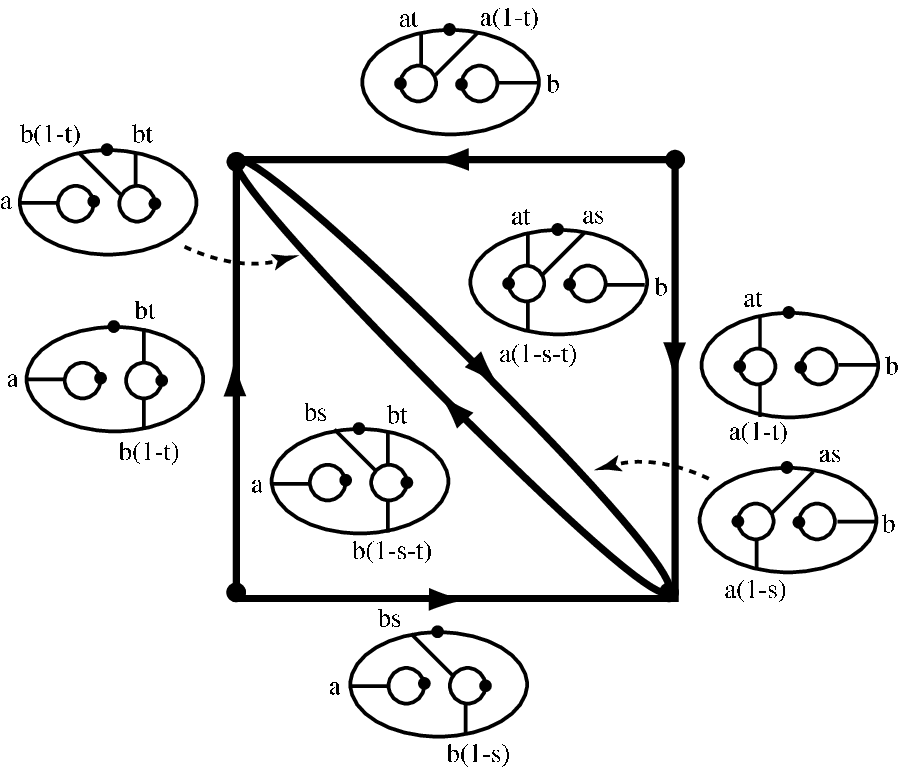}}

\vskip .2in

\centerline{
\begin{footnotesize}
{{\bf Figure 9} Closed sector relations: compatibility of bracket and
composition.}
\end{footnotesize}}

\centerline{\epsffile{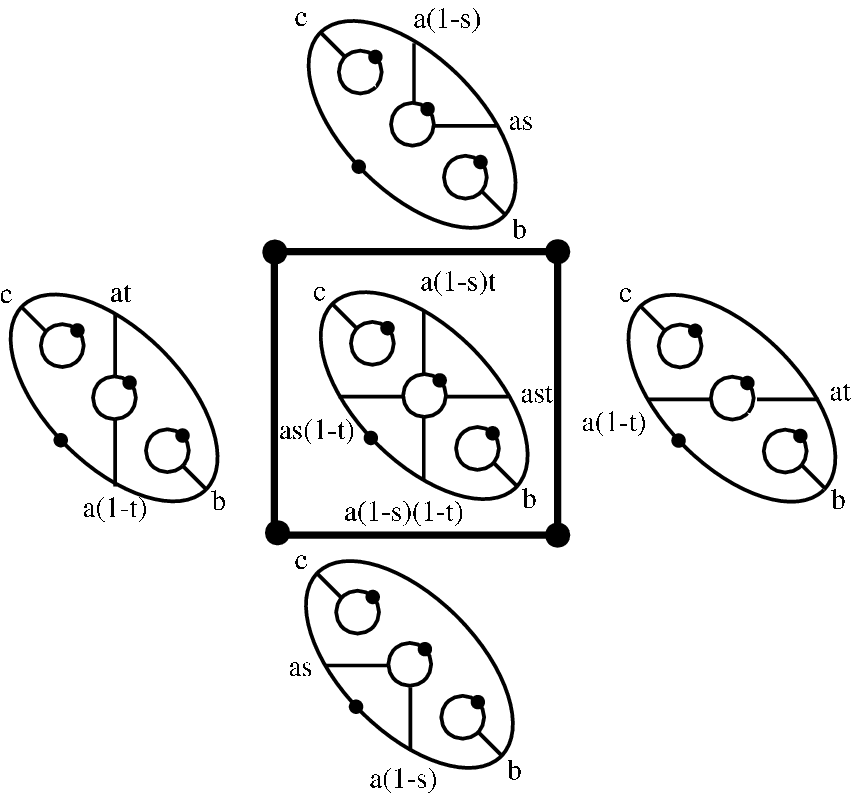}}

\centerline {\begin{footnotesize}{\bf Figure 10} Homotopy for one-third of the BV equation.\end{footnotesize}}

\vskip .2in

The other closed sector relations were already confirmed in \cite{KLP} and are rendered in Figures~8-11
in the current formalism, where all brane-labelings are tacitly taken to be $\emptyset$; furthermore all
boundary-parallel pants curves are omitted from the figures (except in Figure~11 for clarity).
The Frobenius equation is
again degree one, and the BV equation itself is degree two.

\centerline{\epsffile{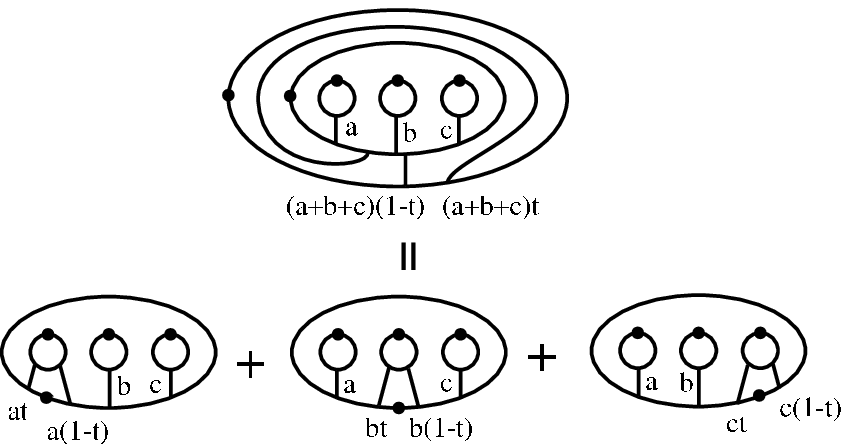}}

\vskip .2in

\begin{footnotesize}

\noindent{\bf Figure 11} The homotopy BV equation.
There are three summands in Figure~11, each of which is parameterized by an
interval, which together combine to give the sides of a triangle.
For each side of this triangle, there is the homotopy depicted in Figure 10 with the appropriate
labeling.  Glue the three rectangles
from Figure~10 to the three sides of
the triangle in the manner indicated. The BV equation is then the fact that
one boundary component of this figure (9 terms) is homotopic to the other
boundary component (3 terms); Figure~11 of \cite{KLP} renders this entire homotopy.

\end{footnotesize}

\vskip .2in

\centerline{\epsffile{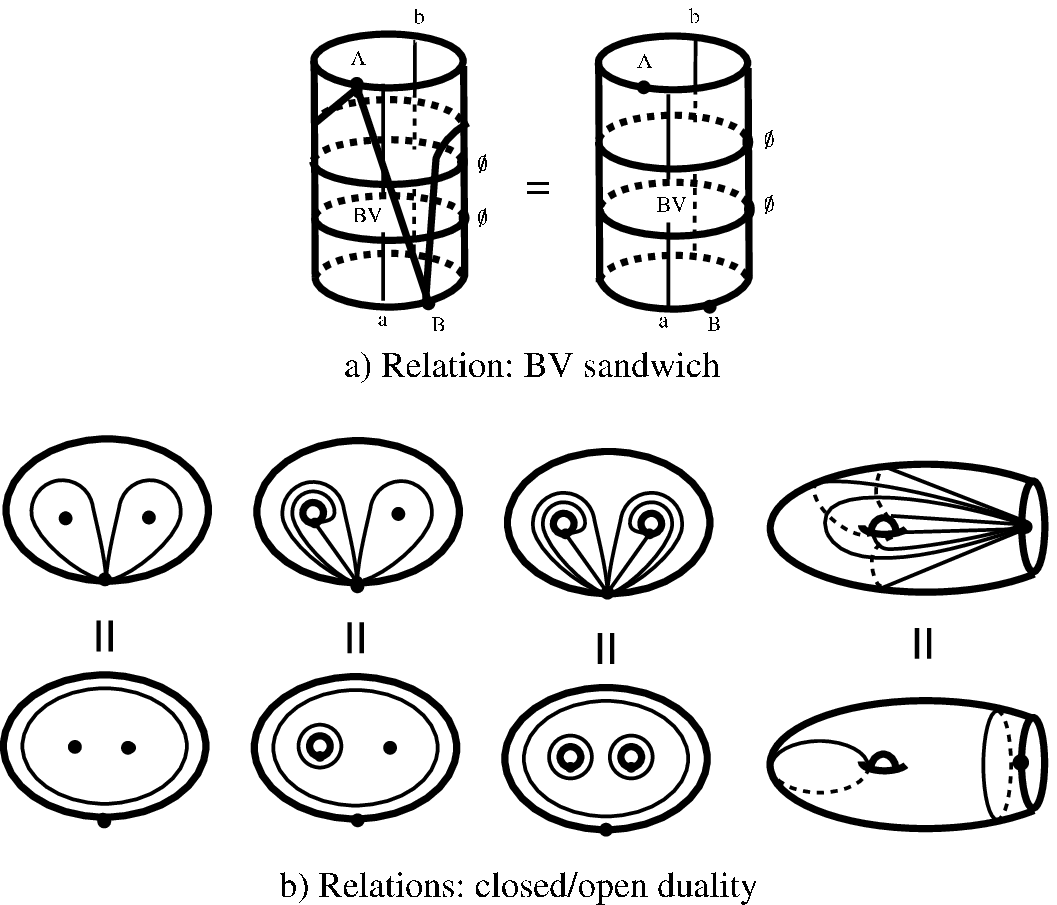}}

\vskip .2in

\centerline{\begin{footnotesize}{\bf Figure 12}~New relations.\end{footnotesize}}

\vskip .2in

It is thus straight-forward to discover new relations, and several
such relations of some significance are indicated in Figure~12.
Figure~12a illustrates a degree one equation on $F^0_{0,2}$ called
the  ``BV sandwich'', which can be  succinctly described by
``close an open string, perform a BV twist, and then open the
closed string''. One-parameter families of weighted arc families
on two  triangles are combined by parameterized open string gluing
to produce a closed string BV twist sandwiched between
closing/opening the string.  The significance of the relations in
Figure~12b and the justification for the choice of terminology
will be explained in the next section.

Here is an algorithm for deriving all of the relations in degree zero on the topological level: Induct over the
topological type of the surface and over the $MC(F)$-orbits of all pairs of generalized pants decompositions of it.
(Though there are only finitely many $MC(F)$-orbits of singleton generalized pants decompositions, there are
infinitely
many $MC(F)$-orbits of pairs.)  In each indecomposable piece, consider each of the possible building blocks
illustrated
in Figures~3-4.  Among  these countably many equations are all of the degree zero equations of the topological c/o
structure.

To derive all higher degree relations on the topological level,
notice that each indecomposable surface has (the geometric
realization of) its arc complexes of some fixed rather modest
dimension.  Thus, parameterized families may be described as
specific parameterized families in each building block, for
instance, in the coordinates of Corollary ~\ref{DTcor}. Such
parameterized families can be manipulated using known
transformations (see the next section) to explicitly relate
coordinates for different generalized pants decompositions and
derive all topological relations.

{\sl A fortiori}, topological relations hold on the chain level (and
likewise for the chain and homology levels as well).  For
parameterized families, there is again the analogous exhaustively
enumerative algorithm, but one must recognize when two
parameterized families are homotopic, which is another level of
complexity.

It is thus not such a great challenge to discover new relations in this manner.  The remaining difficulties involve
systematically understanding not only higher degree equations like the BV sandwich but also in determining a minimal
set
of relations, and especially in understanding the descent to homology.

\vskip .2in

\subsection{Open/closed duality}

\vskip .2in

We seek a collection of combinatorially defined transformations or
``moves'' on generalized pants decompositions of a fixed
brane-labeled windowed surface, so that finite compositions of
these moves act transitively.  In particular, then any closed
string interaction (a standard pants decomposition of a windowed
surface brane-labeled by the emptyset) can be opened with the
``opening operator'' $(i_a^A)^*$ illustrated in Figure~6d, say with a
single brane-label $A$; this surface can be quasi triangulated,
giving thereby an equivalent description as an open string
interaction.

In particular, the two
moves In Figure~13a-b were shown in \cite{P3} to act transitively
on the quasi triangulations of a fixed surface, and likewise
the two ``elementary moves'' on $F^0_{1,1}$ and
$F^0_{0,4}$ of Figure~13c-d were shown in \cite{HT}
to act transitively on standard pants decompositions, where we include also the generalized
versions of Figure~13d on $F_{0,r}^s$ with $r+s=4$ and Figure~13c on $F_{1,0}^1$ as well
(though this includes some non-windowed surfaces strictly speaking).

\vskip .2in

\centerline{\epsffile{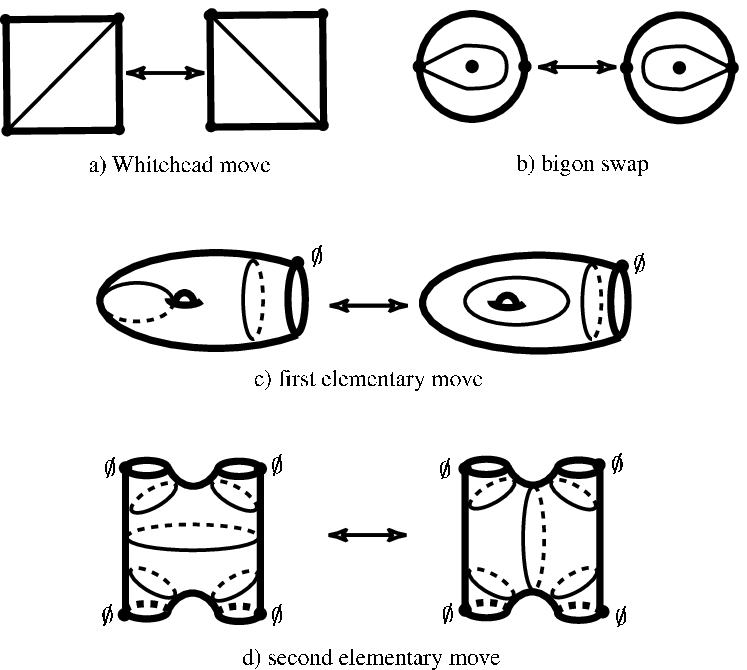}}

\vskip .2in

\centerline{\begin{footnotesize}{\bf Figure 13}~Four combinatorial
moves, where absent brane labels are arbitrary.\end{footnotesize}}

\vskip .2in

For another example, the Cardy equation
can be thought of as a move between the two generalized
pants decompositions depicted in Figure~7e, and likewise for the
four new relations in Figure~12b.

\begin{thm} \label{thmtrans} Consider the following
set of combinatorial moves: those illustrated in Figure~13 together
with the Cardy equation Figure~7e, and the four closed/open duality relations
Figure~12b. Finite
compositions of these moves act transitively on the set of all
generalized pants decompositions of any surface.
\end{thm}

\begin{proof}
In light of the transitivity results mentioned above by topological induction, it remains only to show that the
indicated
moves allow one to pass between some standard pants decomposition and some quasi triangulation of a fixed surface $F$
of type $F_{g,r}^s$.
This follows from the fact that on any surface other than those in Figure 7e and 12b, one can find in
$F$ a separating curve
$\gamma$ separating off one of these surfaces.  Furthermore, one can complete
$\gamma$ to a standard pants decomposition $\Pi$ so that there is at least one window in the same component of
$F-\cup\Pi$ as $\gamma$.  Choose an arc $a$ in $F-\cup\Pi$ connecting a window to $\gamma$; the boundary of a regular
neighborhood of $a\cup\gamma$ corresponds to one of the enumerated moves, and
the theorem follows by induction.
\end{proof}

\vskip .2in

It is an exercise to calculate the effect of these moves on the natural coordinates in Theorem \ref{DT} in
the case of the quadrilateral $F_{0,(4)}^0$ and the once-punctured monogon $F_{0,1}^1$, and the calculation of the new
duality relations on generalized pairs of pants $F_{0,r}^s$, for $r+s=3$, and on $F_{1,1}^0$ is implicit in Figure~4.
The calculation of the first elementary move on $F_{1,1}^0$ is also not so hard, but the formulas are unfortunately
incorrectely rendered in \cite{PH}; see \cite{Pthesis} or \cite{Ptsurvey}. 
The calculation of the second elementary
move on $F_{0,4}^0$, a problem going back to Dehn, was solved in
\cite{Pthesis}.

\vskip .3in

\section{Algebraic properties on the chain and homology levels}

\vskip .2in

\subsection{Operations on the chain level}\label{chainops}

The moves discussed in the last section give rise to relations on
the chain level as well. As explained in Appendix A upon fixing a
chain functor $\Ch$, a chain may be thought of as a parameterized
family of arc families, i.e., as a suitable continuous function
$a(s)\in \A(n,m) $, where $s$ represents a tuple of parameters.
The gluing operations on the chain level can furthermore be
thought of as gluing in families, where the gluing is possible
when the weights of the appropriate windows agree in the two
families. As mentioned previously, any relation on the topological
level gives rise to a relation of degree zero on the chain level.
Some of the relations we discuss will be only up to homotopy,
i.e., of higher degree.

Given any $a\in \Ch(\A(n,m)) $, i.e., given any suitable
parameterized family $a(s)\in\A(n,n)$ on a component of $\A(n,m)$,
say with underlying surface $F$, we may fix a window $w$ on $F$
and regard $a,w$ as an operation in many different ways:

\vskip .1in

\noindent ($n+m-1$)-ary operation:
given chains $a_1,\dots, a_{n+m-1}$, each on a surface with distinguished window, glue them to
all windows of $F$ except $w$; the $n+m-1$ inputs $a_1,\dots, a_{n+m-1}$ yield the
output $w$;

\vskip .1in

\noindent dual unary operation:
given a chain $b$ on a surface with distinguished window, glue it to $F$ along $w$ to produce
the chain we shall denote $a^*(b)$; the input $b$ yields the output $a^*(b)$.

\vskip .1in

More generally, we may partition the windows of the underlying
surfaces into inputs and outputs to obtain more exotic operations
associated with chains.  (The mathematical structure of PROPs were
invented to formalize this structure; for a review see
\cite{MSS,Markl}.)

For instance, let us explain the sense in which the constant chain
$m_{ab}^{ABC}$
in Figure~6b describes the binary operation of multiplication.
Taking the base of the triangle as the distinguished window $w$,
consider families $a(s)\in \Ch(\A(n,m))$ with distinguished window $w_a$ and $b(t)\in \Ch(\A(n',m'))$ with
distinguished window $w_b$, where the brane labels at the endpoints of $w_a$ are $A,B$ and of $w_b$ are $B,C$.
These chains can be
glued to the constant family $m_{ab}^{ABC}$ if and only if
the weight of $a(s)$ on its distinguished window is constant equal to
$a$ and the weight of $b(t)$ on its distinguished window is constant
equal to $b$.
Let the base of
the triangle in Figure 6b be the window 0, the side $AB$ the
window 1, and the side $BC$ the window 2. The  chain operation is
defined as $(m_{ab}^{ABC}\bullet_{1,w_a}a(s))\bullet_{2,w_b}b(t)$.
Notice that the resulting chain will have constant weight $a+b$
on its window 0.  It is in this
sense that we shall regard the constant chain $m^{ABC}_{ab}$ as a binary multiplication.

On the other hand $m^{ABC}_{ab}$ acts as a co-multiplication as well:
given $a(s)$ with brane labels $A,C$ on its distinguished window, we have $m^{ABC*}_{ab}\circ
a(s)=m^{ABC}_{ab}\bullet_{0,w_a}a(s)$.

\subsubsection{Degree 0 indecomposables and relations}
Degree zero chains are generated by zero-dimensional
families, that is, by points of the spaces $\A(n,m)$.

For the
indecomposable brane-labeled surfaces of \S \ref{indec},
the relevant degree 0 chains are enumerated in Figure 6a-e and
6l-n. They become explicit operators by fixing the distinguished window $w$ to be the
lower side in 6a, the base in 6b and the outside boundary in 6c-e
and 6l-n: this is the algebraic meaning of the illustrations in Figure 6.
For
example, 6b and 6e give the respective open and closed binary multiplications $m_{ab}^{ABC}$
and $m_{ab}$, while 6a and 6c give the respective identities
$id_a^{AB}$ and $id_a$ on their
domains of definition, namely, families whose weight on the distinguished window is
constant equal to $a$.
The subscripts indicate compatibility for chain gluing and self-gluing
in the chain level c/o structure, and the superscripts denote brane-labels in the open sector.

It follows from
Figures~7a and 8a, that the multiplications $m^{ABC}_{ab}$ and
$m_{ab}$ are associative:
\begin{eqnarray}
m_{a+c,b}^{ADB}\circ (m_{ac}^{ACD}\otimes
id^{BC}_b)&=&m_{a,b+c}^{ACB}\circ
(id_a^{AC}\otimes  m_{c,b}^{CDB}),\\
m_{a+c,b}\circ
(m_{ac}\otimes id_b)&=&m_{a,b+c}\circ (id_a\otimes m_{c,b}),
\end{eqnarray}
where $\circ$ means
the usual composition of operations.
The dual unary operations to these multiplications satisfy the Frobenius
equations up to homotopy as shown in Figure~7c for the open sector:

\begin{multline}
(id^{AD}_{a} \otimes m^{DBC}_{bc}) \circ (m^{BDA*}_{ab} \otimes id^{BC}_{c})=m^{CDA*}_{a,b+c}
\circ m^{ABC}_{a+b,c} \sim m_{ac}^{CDA*}\circ m^{ABC}_{ac}\\
\sim m^{CDA*}_{a+b,c}
\circ m^{ABC}_{a,b+c}=(m^{ABD}_{ab}\otimes id^{CD}_{c})\circ (id^{AB}_{a}\otimes m^{CDB*}_{bc}),
\end{multline}

and in Figure~8b for the closed sector:
\begin{multline}
(id \otimes m_{bc}) \circ (m^{*}_{ab} \otimes id)=m^*_{a,b+c}
\circ m_{a+b,c} \sim m_{ac}^*\circ m_{ac}\\
\sim m^*_{a+b,c}
\circ m_{a,b+c}=(m_{ab}\otimes id)\circ (id\otimes m_{bc}^*).
\end{multline}

The ``closing'' operation $i_a^A$ of Figure 7d acts as a unary operation which
changes one window from open to closed. Its dual ``opening'' operation 
changes one closed to one open window. It follows from Figure 7b
that $i_a^A$ is an algebra homomorphism:
\begin{equation}
i^A_{a+b}\circ m_{ab}= m^{AAA}_{ab}\circ (i^A_a, i^A_b).
\end{equation}
The image of $i_a^A$ lies in
the center, as in Figure 7d:
\begin{equation}\label{ourcenter}
m^{ABA*}_{ab}\circ i^A_a=\t_{1,2}\circ m^{BAB*}_{ba}\circ i^B_b,
\end{equation}
where $\tau_{12}$ interchanges the tensor factors, namely, interchanges the two non-base sides of the
triangle, and it satisfies the Cardy equation
in Figure 7e:
\begin{equation}
i^A_{a+b}\circ i^{B*}_{a+b}= m_{ab}^{ABA} \circ \tau_{1,2} \circ  m^{BAB*}_{ab}.
\end{equation}

The operators in Figures~6l-n are puncture operators, which are
``shift operators for the puncture grading''.

\vskip .2in

\centerline{\epsffile{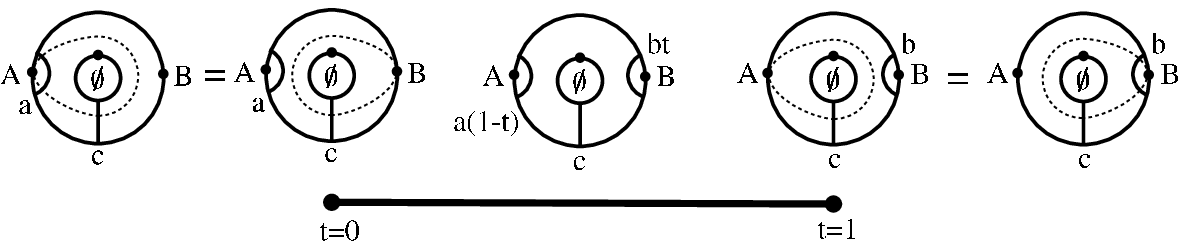}}

\begin{footnotesize}

\noindent{{\bf Figure 14}~Homotopy of Equation (\ref{usualcenter}) to Equation (\ref{ourcenter}), where the dotted
lines indicate generalized pants decompositions.}
\end{footnotesize}

\vskip .2in

We finally express the center equation in a less symmetric but more familiar form:
 \begin{equation}\label{usualcenter}
m_{ac}^{AAB}\circ (i_b^A\circ id)\sim m_{ba}^{ABB}\circ \t_{1,2}
\circ (i^B \circ id).
\end{equation}
\noindent As indicated in Figure~14, the equation (\ref{usualcenter}), which is represented by the
far left and right figures, is equivalent on
the chain level to two copies of our equation (\ref{ourcenter}), represented by the two equalities in the figure,
which holds on the nose.

\vskip .2in

\subsubsection{Degree one indecomposables and relations}\label{degreeoneindecomposables}
The known degree one operations are  given by the binary
operation $\circ_{ab}$ of Figure 7g and the operation $\Delta_a$ of
Figure 7h. These operations are related, indeed, they satisfy the GBV equations up to
homotopy: $\circ_{ab}$ is a homotopy pre-Lie
operation  whose induced homotopy Gerstenhaber structure
coincides with the induced homotopy Gerstenhaber structure of the
homotopy BV operator $\Delta_{a}$ (see Figures 9-11).  This completes the discussion
of know relations.

There is a degree one chain which is of interest, namely, the family which is
generated by the BV sandwich setting $a=t$ and $b=a-t$ in Figure 12a.  This is supported
on the annulus $F_{0,2}^0$ with brane labels A and B. We shall call this operator
$D^{AB}_a$.  The BV sandwich equation then gives the equality of chain level operators
 \begin{equation}\label{firstD}
D^{A,B}_a=(i_a^{B})^*\circ \BV_a \circ i_a^A
\end{equation}

In the same spirit, there is another degree one chain which is
associated to $m_{t,(a-t)}^{ABC}$ for $0\leq t\leq 1$. Although this chain is
not closed, it appears naturally as follows.
Using the BV sandwich relation for the
chain $D_a^{AB}$, we see that it also decomposes as:
\begin{equation}\label{secondD}
D_a^{AB}=m^{ABA}_{(a-t),t}\circ \t_{1,2} \circ m^{BAB*}_{(a-t),t}
\end{equation}

Thus, $D^{A,B}_a$ admits the two expressions (\ref{firstD}) and (\ref{secondD}), so the chain
$D_a^{BC}~D_a^{AB}$, which is a kind of ``BV-squared in the open sector'',
likewise admits the two expressions
corresponding to the two different generalized
pants decompositions of $F_{0,2}^1$:

\begin{eqnarray*}
D_a^{BC}D_a^{AB}&=&(i_a^{C})^*\circ D_a \circ
i_a^B \circ (i_a^{B})^*\circ D_a \circ i_a^A\nn\\
&=&m^{BCB}_{(a-t),t}\circ \t_{1,2} \circ m^{CBC*}_{(a-t),t}\circ
m^{ABA}_{(a-t),t}\circ \t_{1,2} \circ m^{BAB*}_{(a-t),t}.
\end{eqnarray*}

\noindent See the closing remarks for a further discussion of this operator $D_a^{AB}$.

\vskip .2in

\begin{lem}\label{indecompcircle}
Suppose that $(F,\beta )$ is an indecomposable brane-labeled
windowed surface. If $F$ is a triangle or a once-punctured
monogon, then $\widetilde{Arc}(F,\beta )$ is contractible.  For an
annulus, $\widetilde{Arc}(F,\beta )$ is homotopy equivalent to a
circle, and for a generalized pair of pants with $r>0$ boundary
components, $\widetilde{Arc}(F,\beta )$ is homotopy equivalent to
the Cartesian product of $r$ circles.
\end{lem}

\vskip .2in

\begin{proof}  The claims for triangles and once-punctured monogons are clear from Figure~4a-b.
For the degree one
indecomposables, we first have the annuli $F_{0,2}^0$
brane-labeled by $\emptyset,\emptyset$ or by
$\emptyset,A\neq\emptyset$; the free generator of the first homology of
the former is precisely the BV operator $\Delta$, while the
free generator of the latter is $i^A\circ_{1,1} \BV$, where $1$ is the window
labeled by $\emptyset$.
In each case, we have that
$\A(F_{0,2}^0,\beta)$ is homotopy equivalent to a circle.
For $(F_{0,1}^2,\beta_{\emptyset})$, we again have
$\A(F_{01}^2,\beta_{\emptyset})$ homotopy equivalent to a circle as
in Figure~2, with the free generator
$\punc\circ \Delta$.

For the generalized pairs of pants, first notice that the set of {\it all} homotopy classes of families of
projective weighted arcs in  a generalized pair of pants with $r>0$ boundary components (where the arc family need
{\sl not} meet each boundary component) is homeomorphic to the join of $r$ circles.   (In effect, a point in the
circle determines a projective foliation of the annulus as in Figure~3, and one deprojectivizes and combines as in
Figure~4d to produce a foliation of the pair of pants.).  The complement of two spaces in their join is homeomorphic to
the Cartesian product of the two spaces with an open interval, and the lemma follows.
In fact,
the first homology of
$\A(F_{0,2}^1,\beta_{\emptyset})$ is freely generated by $\BV\circ_{1,1} P_2$, and $\BV\circ_{1,2}P_2$, and the first
homology of
$\A(F_{0,3}^0,\beta_\emptyset)$ is freely generated by
$\BV\circ_{1,1} m_a$, $\BV\circ_{1,2}m_a$, and
$\BV\circ_{1,3}m_a$.
\end{proof}

\vskip .2in

For a final chain calculation, consider the degree two chain
 defined by $\BV_{sq}^B=
\BV_a\circ  (i_a^{B})^*\circ i_a^B \circ \BV_a,$ which is another
type of ``BV-squared operator in the open sector'' arising on the
surface $F_{0,2}^1$ with brane-labeling $\beta$ given by $\emptyset$ on the boundary
and by $B$ at the puncture. In fact, $\BV_{sq}^B$ generates
$H_2(\A(F^1_{0,2},\beta))$, where
$H_2(\A(F^1_{0,2},\beta)={\mathbb Z}$ by
Lemma~\ref{indecompcircle}.

Tautologically, $\BV_{sq}^B(s,t)$ can be written as the sum of two
non-closed chains $\BV_{sq}^B=(\BV_{sq}^B)_+ + (\BV_{sq}^B)_-$ given by
\begin{eqnarray}
(\BV_{sq}^B)_+=&\BV _a(t) \circ i_a^B \circ (i_a^{B})^*\circ
\BV_a(s); \quad s+t\leq 1,\nn\\
(\BV_{sq}^B)_-=&\BV _a(t) \circ i_a^B \circ (i_a^{B})^*\circ
\BV_a(s); \quad s+t\geq 1.
\end{eqnarray}

Furthermore, we may homotope each of the operators
$(i_a^{C})^*\circ (\BV_{sq}^B)_+\circ i_a^A$ and $(i_a^{C})^*\circ
(\BV_{sq}^B)_-\circ i_a^A$ into ``traces over multiplications'' in the following
sense, where we concentrate on $(\BV_{sq}^B)_+$ with the parallel discussion
for $(\BV_{sq}^B)_-$ omitted.
Consider the homotopy of arc families in $F_{0,2}^1$
depicted in Figure 15, which begins with $(\BV_{sq}^B)_+$ and ends with the indicated
family.
Cutting on the dotted lines in Figure~15 decomposes each surface into
a hexagon, and these hexagons may be triangulated into four triangles corresponding to four multiplications.
Thus, each of the operations $(\BV_{sq}^B)_+$ and $(\BV_{sq}^B)_-$
is given as the double trace over a quadruple multiplication.
Again, see the closing remarks for a further discussion of these operators.

\vskip .2in

\centerline{\epsffile{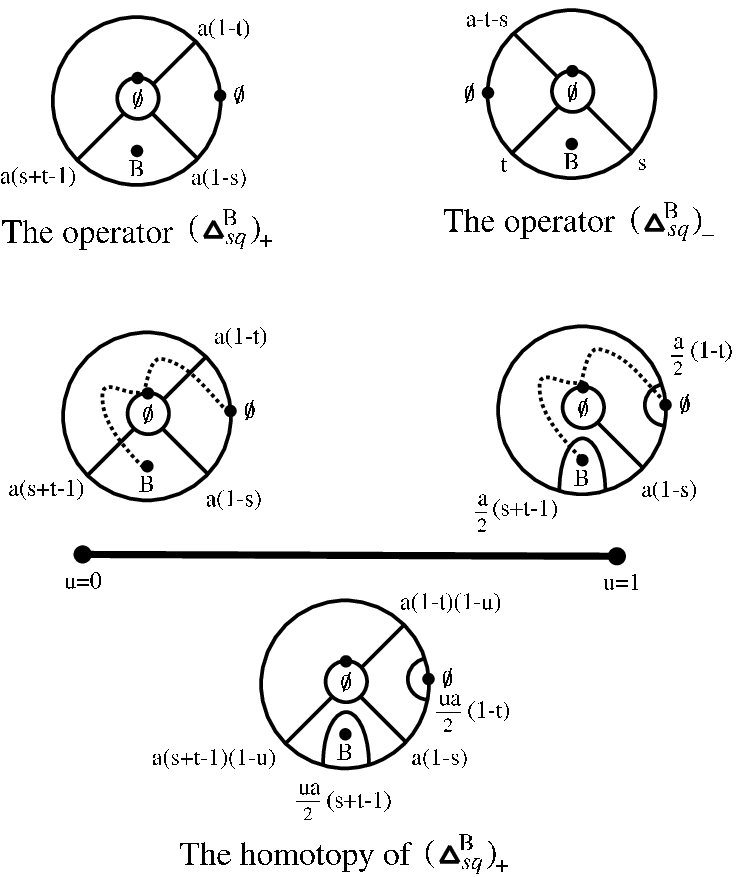}}

\vskip .2in

\begin{footnotesize}
\centerline{{\bf Figure 15}~The operators $(\BV_{sq}^B)_\pm$ and the homotopy of $(\BV_{sq}^B)_+$.}
\end{footnotesize}

\vskip .2in

\subsection{Algebraic properties on the homology level}

Since $\A(F,\beta )$ is connected for any windowed surface $F$ with brane-labeling $\beta$, we conclude
$$H_0(\A(n,m))=\bigoplus {\mathbb Z},$$
the sum over all homeomorphism classes of brane-labeled surfaces $(F,\beta )$ with $n$ closed and $m$ open
windows. It
follows that the degree zero relations on the homology level are precisely those holding on the chain level
up to homotopy.

This observation together with Theorem \ref{thmtrans} implies the
result of \cite{Lew,LP} that the open/closed cobordism group admits the
standard generators with the complete set of relations depicted in
Figures~7 and 8:  associativity, algebra homomorphism, the
Frobenius equations, center and Cardy together with duality.

\begin{lem}
\label{contract} Each component of $\A(0,m)$ is contractible, hence
the homology $H_*(\A(0,m))$ is concentrated in degree $0$.
\end{lem}

\begin{proof}
Consider the foliation which has a little arc around each of the
points of $\delta$ with constant weight one. We can define a flow
on $\A(0,m)$, by including these arcs with any element $\a\in
\A(0,m)$ and then increasing their weights to one while
decreasing the weights of all the original arcs to zero.
\end{proof}

\vskip .2in

\noindent In particular, the first open BV operator $D_a^{AB}$
itself thus vanishes on the homology level, while the second open
BV operator $i_B\circ\Delta$ and even its ``square'' $\Delta
_{sq}^B$ do not.  (The situation may be different in conformal
field theory as discussed in the closing remarks.)

Let $m,m^{ABC},i^A$ be the respective images in homology of the chains
$m_1,m_{11}^{ABC},i^A_1$,  define $\pa,\punc,\pp$ to be the
images of the puncture operators $\pa_1,P^1_{2\; 1},P^2_{1~11}$ in homology, and $\BV$ the
image in homology of $\Delta_1$.

Just as chains can be regarded as operators on the chain level, so too homology classes can be regarded
as operators on the homology level.

\begin{prop}
The degree zero operators on homology are precisely generated by the degree zero
indecomposables $m,m^{ABC}$ and $i^A$
provided $\emptyset\not\in\beta (\sigma )$, where $\sigma$ denotes the set of punctures.
If $\emptyset\in\beta (\sigma )$, then one must furthermore include the operators
$\pa,\punc,\pp$.  The degree zero relations
on homology are precisely those given by the moves of Theorem
\ref{thmtrans}.
All operations of all degrees
supported on indecomposable surfaces are
generated by the degree zero operators and $\Delta$.
\end{prop}

\begin{proof}
The degree zero operators arise from constant families and each $\A(F,\beta )$ is connected.  By
Theorem
\ref{DT}, the operators of degree zero arise from degree zero chains
on indecomposable surfaces, proving the first part.
The second part follows from Lemma~\ref{indecompcircle}.
\end{proof}

\begin{thm}\label{homthm} Suppose $\emptyset\notin\beta (\sigma )$.  Then an algebra over the modular
bi-operad $H_*(\coprod_{n,m}\A(n,m))$ is a pair of vector spaces
$(C,\Ass)$ which have the following properties: $C$ is a
commutative Frobenius BV algebra $(C,m,m^*,\BV)$, and
$\Ass=\bigoplus_{(A,B\in {\PB\times \PB})}\Ass_{AB}$  is a
$\PB$-colored Frobenius algebra (see e.g., \cite{LP} for the full
list of axioms). In particular, there are multiplications
$m^{ABC}:\Ass_{AB} \otimes \Ass_{BC}\to \Ass_{AC}$ and a
non-degenerate metric on $\Ass$ which makes each $\Ass_{AA}$ into
a Frobenius algebras.

Furthermore, there are morphisms $i^{A}:C\rightarrow \Ass_{AA}$
which satisfy the following equations: letting $i^*$ denote the dual of
$i$, $\t_{12}$ the morphism permuting two tensor factors, and letting
$A,B$ be arbitrary non-empty brane-labels, we have
\begin{eqnarray}
&i^B\circ {i^A}^* =  m_B \circ \t_{12}\circ m^*_A& (\text{Cardy})\\
&i^A(C)~{\rm is~central~in}~{\mathcal A}_A &(\text{Center})\\
& i^A\circ \BV \circ
{i^B}^*=0& (\text{BV vanishing})
\end{eqnarray}
These constitute a spanning set of operators and a complete set of
independent relations in degree zero.  All operations of all degrees
supported on indecomposable surfaces are
generated by the degree zero operators and $\Delta$.
 \end{thm}

\begin{proof}
By definition, an algebra over a modular
operad is a vector space with a non-degenerate bilinear form such
that the operations are compatible with dualization \cite{GK}.
The previous proposition substantiates the first sentence of the theorem, and the second sentence
follows in particular.

The claim that $i^A$ and $i^*_A$ are morphisms
and the Cardy and center equations then follow from the chain level equations
and the fact that the operation $a^*$ for a chain $a$ is
in fact the dual operation on the Frobenius algebra. Indeed by
Theorem \ref{thmtrans}, the Cardy and Center equations generate the relations
in degree zero.
Finally, BV vanishing follows from Lemma~\ref{contract} and the last assertion
from Lemma~\ref{indecompcircle}.
\end{proof}

\vskip .2in

Notice that if $A$ is a symmetric Frobenius algebra with a pairing
$\langle\;,\;\rangle$ then our center equation (\ref{ourcenter}) $m^{ABA*}\circ
i^a=\t_{1,2}\circ m^{BAB*}\circ i^B$ implies that $i^A$
takes values in the center, that is, the equation $m^{AAB}\circ
(i^A\circ id)=m^{ABB}\circ \t_{1,2} \circ (i^B \circ id)$ holds.
Indeed, we have
\begin{multline}
\la i^A(a)b,c\ra=\la i^A(a),bc\ra=\la m^{AAA*}(a),b\otimes c\ra\\
=\la \t_{1,2} \circ m^{AAA*}\circ i^A(a),b \otimes c\ra=\la
a,cb\ra=\la ba,c\ra,
\end{multline}
where the last equation holds since the Frobenius algebra was
assumed to be symmetric.

There is a further grading
by the number of punctures which are brane-labeled by $\emptyset$.
(Recall that gluing operations never give rise to new punctures labeled by $\emptyset$.)
In this case, representations will actually lie in triples of vector spaces 
$(C,\Ass,W)$, where $W$ 
corresponds to closed string
insertions which give deformations of the original operations.

 \vskip .3in

\section{Closing remarks}

An important challenge is to understand the transition to conformal field theory from 
the topological field theory described in this
paper
as the degree zero part of the homology.

The higher degree part of the homology studied in the body of this paper is based on ``exhaustive''Ê 
arc families $\A(F,\beta )$ which meet each window of a brane-labeled windowed surface $(F,\beta )$.
As mentioned in the Introduction, such an exhaustive family isÊ 
not enough to determine a metric on
the surface unless the arc family
quasi fills the surface, i.e., complementary regions are either polygons or exactly once-punctured polygons, and we
shall let $\A_\#(F,\beta )\subseteq \A(F,\beta )$ denote the corresponding subspace.

Furthermore as in the Introduction, the subspace $\A_\#(F,\beta )$ 
is naturally identified with
what is essentially Riemann's moduli space of $F$ with one
point in each boundary component.  
One might thus hope to describe 
CFT via the chains and homology groups
of the quasi filling subspace $\A_\#(n,m)$ of the spaces $\A(n,m)$.  
In
this setting, both the chain and homology levelsÊ 
must be re-examined compared to the exhuastive case studied in the body 
of the paper
insofar as homotopies must respect
$\A_\#(n,m)$, which in particular is not
invariant under the gluing or self-gluing operations.
To give a c/o-structure on the topological or chain 
level on $\A_\#(n,m)$, one can
imagine using homotopies in the appropriate combinatorial
compactification (see \cite{P1}) to define the gluings
for compactified quasi filling arc families, or alternatively, one might
proceed solely on the level of
cellular chain complexes, see below.

By Lemma \ref{contract} for exhaustive arc families, the open sector BV
operator $D^{BB}=i^{B*}\circ \Delta\circ i^{B}$ and its square vanish on theÊ 
level of homology, while
on the other hand, the other open BV-squared operator $\Delta ^B_{sq}=\Delta
\circ i^B\circ i^{B*}\circ \Delta$ is not zero, but rather aÊ 
generator of
the second homology group of $\A(F_{0,2}^1,\beta)$, where
$\beta$ takes value $\emptyset$ on the boundary and value $B$ at theÊ 
puncture.  In physical terms for exhaustive families, the corresponding bulk operator $\Delta ^B_{sq}$ vanishes only
after coupling to the boundary, i.e., $(D^{BB})^2=0$ yet $\Delta ^B_{sq}\neq 0$.

On the other hand in the context of quasi filling arc families,
Lemma \ref{contract} does {\sl not} hold, and neither the operator
$D^{BB}$ nor its square now vanishes.
This serves to emphasize one basic algebraic
difference between exhaustive and quasi filling arc families,
which are presumably required for CFT.

Furthermore, the
non-vanishing of $\Delta _{sq}^B$
is reminiscent of the
appearance of the Warner term in Landau-Ginzburg theory \cite{Warner}.
The other open sector BV-squared operator $(D^{BB})^2$ we consider is supported on the
surface $F_{0,2}^1$ that is the ``open square''Êof the surface $F_{0,2}^0$ which supports
the Cardy equation.  This gives additional credence to this point of viewÊ 
since it is shown in
\cite{Bru1} that the Cardy condition is intimately related to the
compensatory term required to make the action of the LG model BRST
invariant.

The fact that the open sector BV-squared operator $(D^{BB})^2$ vanishes
in the exhaustive case is what one would naively expect.
However, the non-vanishing of this operator in the quasi filling caseÊ 
might help to explain the
appearance of unexpected D-branes in the LG models,
cf. \cite{Bru1,Laz2,KR}. In particular in relation to
Ê Kontsevich's approach to $D$-branes on LG-modes (cf.\cite{KR}),
one might ask if our open sector BV-squared operator $\Delta_{sq}^B$ or
another ``square of a BV-like operator" satisfies an equation of
the form $\Delta_{sq}^B~a=[U_B,a]$, or in a representation,
$\Delta_{sq}^B~m=U_B~m$, for some operator $U_B$, on the level of
either exhaustive or quasi filling families.

These remarks explain
our attention to the chains $D_a^{A,B}$ and $\Delta _{sq}^B$ in
Section~\ref{degreeoneindecomposables}.
Although our results do
not match this formulation exactly, the
decomposition $\Delta _{sq}^B=(\Delta ^B_{sq})_+ + (\Delta ^B_{sq})_-$ is suggestive of a
commutator equation, and the homotopy we described in Section~\ref 
{degreeoneindecomposables}
shows that each of $(\Delta ^B_{sq})_\pm$ is indeed a sort of multiplication operator, forÊ 
instance,
if the representing algebra is
super-commutative.

In summary, the fact that the open sector BV-squared
operator $(D^{BB})^2$ does not vanish in the
quasi filling case and that $\Delta_{sq}^B$ does not vanish in theÊ 
exhuastive case may
be regarded as the statement that the boundary contribution of the
BRST operator need not square to zero; rather, it may be necessary to
introduce additional terms to make the entire action BRST invariant.
Further analysis could give conditions on the
representations of algebras over our c/o structure which can beÊ 
considered
physically relevant.

Although it has not been the focus in this paper, we wish to point out
that there is a discretized version of $G$-colored c/o structures in the
category of topological spaces giving $G$ the discrete
topology, and the $G$-colored c/o-structure then naturally descends to a
$G$-colored c/o structure on both the chain and
homology levels. We can for instance restrict the
topological c/o structure on $\A(F,\beta )$ to the subspaces
where each window has total weight given by a natural number.

In particular for closed strings, it can be shown \cite{hoch}
that there is a natural chain complex of open or relative cells,
which calculates the homology of the moduli spaces, i.e., the homology
of $\A_\#(F,\beta_\emptyset )$ which can be
given the structure of an operad.
In effect, these spaces are
graded by the number of arcs in an arc family,
and the corresponding filtration is preserved by the
gluing operations when viewed as operations on filtered families. Now
projecting to the associated graded object of the filtration, oneÊ 
obtains a cell
level operad.
Furthermore, discretizing as in the previous paragraph, 
one obtains actions on the tensorÊ 
algebra and on
the Hochschild co-chain complex of a
Frobenius algebra.
This discretized and filtered elaboration of the c/o structure could give a
formulation of a version of CFT purely in terms of algebraic topology.
A proving ground for these
considerations might be the topologicalÊ LG models of
\cite{Vafa,Man3c} and their orbifolds \cite{IV,orbsing}.

Let us also mention that Thurston invented a notion of ``tangential measure''
(see \cite{PH}) precisely to capture the lengths as opposed to the
widths of the rectangles in a measured foliation, suggesting yet another geometric aspect of this
passage from TFT to CFT.Ê It is perhaps also worth saying
explicitly that an essential point of Thurston theory is that
twisting about a curve accumulates in the space of {\sl
projective} measured foliations to the curve itself, and this
suggests that the limit of the BV operator $\Delta (t)$ as $t$
diverges might be profitably studied projectively in the context
of Appendix~B.

There is presumably a long way to go until the algebraic structureÊ 
discovered
here is fully understood in higher degrees on the level of homology,Ê 
let alone
for {\sl compactifications} of arc complexes in the quasi filling andÊ 
exhaustive cases.
In the quasi filling case on the level of homology, the underlying
groups supporting these operations comprise the homologyÊ 
groups of Riemann's moduli spaces of bordered surfaces, which are
themselves
famously unknown, yet these unknown groups apparently support the modular bi-operad structure ofÊ 
Theorem~\ref{homthm} at least in this discretized filtered sense.

\vskip .3in

\renewcommand{\theequation}{A-\arabic{equation}}
\renewcommand{\thesection}{A}
\setcounter{equation}{0}  
\setcounter{subsection}{0}

\section*{Appendix A: c/o structures}

\vskip .2in

\subsection{The definition of a c/o-structure}
\label{coapp}
 Specify an object ${\mathcal O}(S,T)$ in some fixed
symmetric monoidal category for each pair $S$ and $T$ of finite
sets.  A {\it $G$--coloring} on ${\mathcal O}(S,T)$ is the further
specification of an object $G$ in this category and a morphism
$\mu : S\sqcup T \to Hom({\mathcal O}(S,T),$G$)$, and we shall let
${\mathcal O}_\mu (S,T)$ denote this pair of data.

A $G$-colored ``closed/open'' or {\it c/o structure} is a
collection of such objects ${\mathcal O}(S,T)$ for each pair of
finite sets $S,T$ together with a choice of weighting $\mu$ for
each object supporting the following four operations which are morphisms in the
category:

\vskip .1in

\noindent
\noindent{\it Closed gluing}:~~$\forall s\in S, \forall s'\in S'$ with $\mu (s)=\mu '(s')$,
$$\circ_{s,s'}: {\mathcal O}_\mu (S,T)\otimes{\mathcal O}_{\mu '} (S',T')\to {\mathcal O}_{ \mu ''}
(S\sqcup S'-\{ s,s'\}
, T\sqcup T');$$

\vskip .1in

\noindent{\it Closed self-gluing}:~~$\forall s,s'\in S$ with $\mu (s)=\mu (s')$ and $s\neq s'$,
$$\circ^{s,s'}: {\mathcal O}_\mu (S,T)\to {\mathcal O}_{\mu ''} (S-\{ s,s'\} , T);$$

\vskip .1in

\noindent{\it Open gluing}:~~
$\forall t\in T, \forall t'\in T'$ with $\mu (t)=\mu '(t')$,
$$\bullet_{t,t'}: {\mathcal O}_\mu (S,T)\otimes{\mathcal O}_{\mu '} (S',T')\to {\mathcal O}_{\mu ''}
(S\sqcup S' , T\sqcup T'-\{t,t'\});$$

\vskip .1in

\noindent{\it Open self-gluing}:~~$\forall t,t'\in T$ with $\mu (t)=\mu (t')$ and $t\neq t'$,
$$\bullet^{t,t'}: {\mathcal O}_\mu (S,T)\to {\mathcal O}_{\mu ''} (S , T-\{ t,t'\}).$$

\noindent In each case, the coloring $\mu ''$ is  induced in the
target in the natural way by restriction, and we assume that $S\sqcup S'\sqcup T\sqcup T'-\{
s,s',t,t'\}\neq\emptyset$.

The  axioms are that the operations are equivariant
for bijections of sets and for bijections of pairs of sets, and the collection of all operations taken together
satisfy associativity.

Notice that we use the formalism of operads indexed by finite
sets rather than by natural numbers as in \cite{MSS} for instance.

\leftskip=0ex

\subsection{Restrictions}
A c/o structure specializes to standard algebraic objects
in the following several ways.

There are the two restrictions
 $(\CO(S,\emptyset),\circ_{s,s'})$ and
 $(\CO(\emptyset,T),\bullet_{\t,\t'})$
each of which forms a $G$-colored cyclic operad in the usual sense.

The spaces
$(\CO(S,T), \circ_{s,s'},\bullet_{\t,\t'})$ with only the non self-gluings as structure maps
form a cyclic
$G\times \Z/2\Z$-colored operad, where the $\Z/2\Z$ accounts for
open and closed, e.g., the windows labeled by $S$ are
regarded as colored by $0$ and the windows labeled by $T$ are
regarded as colored by $1$.

If the underlying category has a
coproduct (e.g., disjoint union for sets and topological
spaces, direct sum for Abelian groups and linear spaces), which we
denote by $\coprod$, then the indexing sets can be regarded as
providing a grading: i.e., \ $(\coprod_T
\CO(S,T),\circ_{s,s'})$ form a cyclic $G$-colored operad graded by
the sets $T$, and
$(\coprod_S\CO(S,T),\bullet_{\t,\t'})$ form a cyclic
$G$-colored operad graded by the sets $S$.

\subsection{Modular properties}
\label{modular}
There is a relationship between c/o structures and modular operads. Recall
that in a modular operad there is an additional grading on the
objects, which is additive for gluing and increases by one for
self-gluing. Imposing this type of grading here, we
define a {\it $(g,\chi -1)$ c/o-structure} to be a c/o structure with
two gradings $(g,\chi)$,
$$\CO(S,T)=\coprod_{g\geq 0,\chi
\leq 0} \CO(S,T;g,\chi)$$ such that

\vskip .1in

\leftskip .2in

\noindent (1)
$\CO(S,T;\chi-1)=\coprod_{g\geq 0}\CO(S,T;g,\chi)$ is additive in
$\chi-1$ for $\bullet_{t,t'}$, and $\chi-1$ increases by one for
$\bullet^{t,t'}$; and

\noindent (2) $\CO(S,T;g)=\coprod_{\chi\leq
0}\CO(S,T;g,\chi)$ is additive in $g$ for $\circ_{s,s'}$, and $g$
increases by one for $\circ^{s,s'}$.

\vskip .1in

\leftskip=0ex

It follows that a $(g,\chi -1)$
c/o structure is a modular $G$-colored bi-operad in the sense that the
 $\CO(S,T;g)$ form a $T$--graded $\Rp$-colored
 modular operad\footnote{We impose
neither  $3g-3+|S|>0$ nor $3(-\chi+1)+|T|-3\geq 0$.} for the  gluings $\circ_{s,s'}$ and
$\circ^{s,s'}$, and the $\CO(S,T;1-\chi)$ form an $S$--graded
$\Rp$-colored  modular operad$^*$ for the gluings $\bullet_{t,t'}$ and
$\bullet^{t,t'}$.

\subsection{Topological and Chain level c/o structures} A
 {\it topological c/o structure} is an
$\Rp$-colored c/o structure  in the category of topological spaces.

Let $\Ch$  denote a chain functor  together with fixed functorial
morphisms for products $P:\Ch(X)\otimes \Ch(Y)\rightarrow
\Ch(X\times Y)$, {\sl viz.} a chain functor of monoidal
categories. The chain functors of cubical or PL chains for
instance come naturally equipped with such maps, and for
definiteness, let us just fix attention on PL chains. A {\it chain
level c/o structure} is a $\Ch(\Rp)$-colored c/o-structure in the
category of chain complexes of Abelian groups. Notice that if a
collection $\{\CO(S,T)\}$ forms a topological c/o structure, then
we have natural maps
$$\Ch(\mu):S\sqcup T\to Hom(\Ch(\CO(S,T)), \Ch(\Rp)).$$
These maps
together with the induced operations make the collection $\{\Ch(\CO(S,T))\}$ into
a {\it chain level c/o structure} by definition. The compatibility equation for
self-gluings explicitly reads $\Ch(\mu)(w)(a)=\Ch(\mu)(w')(a)$ and
for non self-gluings, we have
\begin{multline}
\label{chaincompat} P(\Ch(\mu)\otimes \Ch(\mu))(w)(P(a\otimes
b))\\=P(\Ch(\mu)\otimes \Ch(\mu))(w')P(a\otimes b)).
\end{multline}
It is not true that a topological c/o structure begets a
$\Rp$-colored structure on the chain level, since the topology on
$\Rp$ is not the discrete topology. For PL chains, we may regard a
generator $a\in \Ch(\CO(S,T))$ as a parameterized family, say
depending on parameters $s$, and we shall denote such a
parameterized family $a(s)$.  The equation (\ref{chaincompat})
then simply reads $\forall s,t$, we have
$\mu(w)(a(s))=\mu(w')(b(t))$.

\subsection{The homology level} The coloring in a topological c/o structure is
given by the contractible group $\Rp$, and we take the
coloring or grading by $\Ch(\Rp)$ in the definition of chain level c/o structure.

On the homology level, the grading
$H_*(\mu):H_*(\CO(S,T))\to H_*(\Rp)$ becomes trivial.  Furthermore, it
is in general not possible to push a c/o structure down to the level of homology
since the gluing and self-gluing operations on the topological or chain level are
defined only if certain restrictions are met. It is, however, possible in
special cases to define operations by lifting to the chain level.

\begin{lem}
\label{scaling} Let  $\CO(S,T)$ be a topological c/o structure
such that each $\CO(S,T)$ is equipped with a continuous $\Rp$ action $\rho$
that diagonally acts on the $\Rp$ grading
$\mu(\rho(r)(\a))=r\mu(\a)$.  Then the homology groups
$H_*(\CO(S,T)$ form a cyclic two colored operad under non self-gluings induced by
$\circ_{s,s'}$ and $\bullet_{t,t'}$.
\end{lem}

\begin{proof} As in Appendix~C, define a continuous flow $\psi _t^w:\CO(S,T)\to \CO(S,T)$ for each $w\in S\sqcup T$
by
$$\psi _t^w(\alpha )=\rho\bigl (1-t+{t/{\mu (w)}}\bigr )~(\alpha );$$
thus, $\psi _0^w$ is the identity, and $\psi _1^w(\alpha )$ has weight one on $w$.
Given two cohomology classes $[a]\in
H_*(\CO(S,T))$ and $[b]\in H_*(\CO(S',T'))$ represented by chains
$a \in \Ch(\CO(S,T))$ and $b\in \Ch(\CO(S',T'))$ as well as two
elements $(w,w')\in (S\times S')\sqcup ( T\times T')$, we use the
flows $\psi_t^w,\psi _t^{w'}$ to move $a$ and $b$ into a compatible position by a homotopy.
Explicitly, defining $\tilde a_t=\Ch(\psi_t^{w})(a)$ and $\tilde
b_t=\Ch(\psi_t^{w'})(b)$, we have
$$\Ch(\mu)(w)(\tilde a_1) =\Ch(\mu')(w')(\tilde b_1)\equiv 1.$$
The condition (\ref{chaincompat}) is therefore met, and we define
$$[a]\circ_{s,s'}[b]=[ \Ch(\circ_{s,s'})P(\tilde a_1, \tilde
b_1)],$$ and likewise for $[a] \bullet_{t,t'} [b]$. Associativity of the operations
follows as in
Lemma~\ref{clem1}.
\end{proof}

\begin{prop}
\label{homdescent}
 Let  $\CO(S,T)$ be a topological c/o structure satisfying
 the hypotheses of Lemma  \ref{scaling}.  Furthermore, suppose that
for each $\CO(S,T)$ and each choice of $w\in S\sqcup T$, there is a continuous flow
$\phi^w_t:\CO(S,T)\to \CO(S,T)$, for $0\leq t <1$,
such that  $\phi _0^w$ is the identity,  and for any other $w\neq
w'\in S\sqcup T$ with $\mu(w')(\a)\leq\mu(w)(\a)$, there is a time
$t_c=t_c(\alpha,w')$ for which
$\mu(w')(\phi_{t_c}(\a))=\mu(w)(\phi_{t_c}(\a))$, where $t_c(\alpha ,w')$
depends continuously on $\alpha$.
Then the homology groups $H_*(\CO(S,T))$ carry operations induced by
$\circ_{s,s'},\circ^{s,s'}$ and $\bullet_{t,t'},\bullet^{t,t'}$.

Moreover, given parameterized families $a,b,c$ and letting $\odot$
denote either operation $\circ$ or $\bullet$, suppose that
$(a\odot ^{u,v} b)\odot ^{v',w} c$ and $a\odot ^{u,v} (b\odot
^{v',w} c)$ are homotopic, that $(a\odot ^{u,v} b)\odot _{v',w} c$
and $a\odot ^{u,v} (b\odot _{v',w} c)$ are homotopic, and that
$(a\odot _{u,v} b)\odot ^{v',w} c$ and $a\odot _{u,v} (b\odot
^{v',w} c)$ are homotopic. Then the operations on $H_*(\CO(S,T))$
are associative.

Finally, if the $\CO(S,T)$ furthermore form a
topological $(g,\chi -1)$--c/o structure, then the induced structure on
homology is a modular bi-operad in the sense of \S\ref{modular}.
\end{prop}

\begin{proof}
The non-self-gluing operations are already present and associative by Lemma
\ref{scaling}.
For the self-gluings, the descent of the operations to homology is described in analogy to
Lemma~\ref{clem3}, and the associativity on the chain and hence homology levels
of the operations finally follows from the assumed existence of the homotopies.
\end{proof}

\subsection{Brane-labeled c/o structures} A  brane-labeled c/o
structure is a c/o structure $\{\CO(S,T)\}$  together with a fixed
Abelian monoid $\mathcal{P}$ of brane labels and for each $\a\in
\CO(S,T)$ a bijection $N_{\a}:T\rightarrow T$ and a bijection
$(\lambda_{\a},\rho_{\a}): T\rightarrow \mathcal{P} \times
\mathcal{P}$, such that

\begin{itemize}
\item[(1)]
 $\rho(t)=\lambda(N(t))$,
 \item[(2)] if
$N_{\a}(t)\neq t$ and $N_{\a'}(t')\neq t'$
\begin{eqnarray*}N_{\a\bullet_{t,t'}\a'}(N_{\a}^{-1}(t))=N_{\a'}(t'),&&
N_{\a\bullet_{t,t'}\a'}(N_{\a'}^{-1}(t'))=N_{\a}(t)\\
\rho_{\a\bullet_{t,t'}\a'}(N^{-1}(t))=\lambda_{\a}(t)\rho_{\a'}(t')&&
\lambda_{\a\bullet_{t,t'}\a'}(N(t))=\lambda_{\a'}(t')\rho_{\a}(t)
\end{eqnarray*}

\item[(3)] $N_{\a}(t)\neq t$ and $N_{\a'}(t')\neq t'$
\begin{eqnarray*}
N_{\bullet^{t,t'}(\a)}(N_{\a}^{-1}(t))=N_{\a}(t'),&&
N_{\bullet^{t,t'}(\a)}(N_{\a}^{-1}(t'))=N_{\a}(t)\\
\rho_{\a\bullet^{t,t'}\a'}(N^{-1}(t))=\lambda_{\a}(t)\rho_{\a}(t'),
&&\lambda_{\a\bullet^{t,t'}\a'}(N(t))=\lambda_{\a}(t')\rho_{\a}(t)
\end{eqnarray*}
\item[(4)] If either $N_{\a}(t)=t$ or $N_{\a'}(t')=t'$ but not
both, then in the above formulas, one should substitute
 $N_{\a'}(t')$ for $N_{\a}(t)$ in the first case and inversely in
 the second case. (If both $N_{\a}(t)=t$ and
 $N_{\a'}(t')=t'$, then there is no equation.)
\end{itemize}

This is the axiomatization  of the geometry given by open windows
with endpoints labeled by right ($\rho$) and left ($\lambda $) brane labels, their order and
orientation along the boundary components induced by the
orientation of the surface, and the behaviour of this data under
gluing.

 For a brane-labeled c/o structure and an idempotent submonoid $\B
 \subset \mathcal{P}$ (i.e., for all $b\in \B, b^2=b$),
 one has the {\rm
$\B\times \B$-colored} substructures defined by restricting the
gluings $\bullet_{t,t'}$ and $\bullet^{t,t'}$ to compatible colors
$\lambda(t)=\rho(t')$.

The relevant example for us is $\B$ the set of branes,
$\mathcal{P}=\mathcal{P}(\B)$ its power set with the operation of
union, where $\B\hookrightarrow \mathcal{P}(\B)$ is embedded by
considering $B\in\B$  as the singleton $\{ B\}$.

\subsection{The  c/o-structure on weighted arc families}
In this subsection we give the technical details for the proof of Theorem
\ref{cothm}.

 First set

\begin{eqnarray*} \A(S,T)&=&
\bigl\{(\a,\phi,\psi)|\a\in \A(|S|,|T|),\\
&&\phi: S \isoto \{\text{Closed
windows}\},\\&& \psi: T \isoto \{\text{Open Windows}\}\bigr\}.
\end{eqnarray*}

Define the respective operations $\circ_{s,s'},\circ^{s,s'},\bullet_{t,t'}$ and
$\bullet^{t,t'}$  to be the closed gluing and self-gluing and open gluing and self-gluing
operations defined in \S\ref{oparcs},
where the
$\Rp$-coloring is  the map $\mu$  given by associating the total weight of a weighted arc family
to a window
$w \in S \sqcup T$.

The $(g,\chi)$-grading is given as follows.  For $\a \in \A(S,T)$,
we let $g$ be the genus and let $\chi$
be the Euler characteristic of the underlying surface $F$; if $F$
has punctures $\sigma$, then by definition $\chi (F)=\chi (F\cup\sigma )-\#\sigma$.

Finally the brane-labeling is given by taking $\lambda$ to be the
brane-labeling of the left boundary point and $\rho$ to be that of the right boundary
point of the window.

\vskip .2in

\vskip .2in

\renewcommand{\theequation}{B-\arabic{equation}}
\renewcommand{\thesection}{B}
\setcounter{equation}{0}  
\setcounter{subsection}{0}
\section*{Appendix B: c/o structure on measured foliations}

\subsection*{{\rm B.1}~ Thurston's theory for closed surfaces}

\vskip .2in

Let us specialize for simplicity in this section to a compact surface $F=F_{g,0}^0$ without boundary with $g>1$
in order to
very  briefly describe Thurston's theory of measured foliations; see \cite{FLP} and \cite{PH} for more detail.

A {\it measured foliation} of $F$ is a one-dimensional foliation ${\mathcal F}$ of $F$
whose singularities are topologically equivalent to the standard $p$-pronged singularities of a
holomorphic quadratic differential $z^{p-2}dz^2$, for $p\geq 3$, together with a
transverse measure $\mu$ with no holonomy, i.e., if $t_0,t_1$ are transversals to ${\mathcal F}$
which are homotopic through transversals keeping endpoints on leaves of ${\mathcal F}$, then $\mu (t_0)=\mu (t_1)$;
$\mu$ is furthermore required to be $\sigma$-additive in the sense that if a transversal $a$ is the countable
concatenation of sub-arcs $a_i$ sharing consecutive endpoints, then
$\mu (a)=\sum _i \mu( a_i)$.

Examples arise by fixing a complex structure on $F$ and taking for the foliation ${\mathcal F}$ the horizontal
trajectories of some holomorphic quadratic differential $\phi$ on $F$.  In the neighborhood of a non-singular point
of
$\phi$, there is a local chart $X:U\subseteq F\to{\mathbb R}^2\approx{\mathbb C}$, so that the leaves of
${\mathcal F}$
restricted to $U$ are the horizontal line segments $X^{-1}(y=c)$, where $c$ is a constant.  If the
domains $U_i,U_j$ of
two charts
$X_i,X_j$ intersect, then the transition function $X_i\circ X_j^{-1}$ on $U_i\cap U_j$ is of the
form $(h_ij(x,y),c_{ij}\pm
y)$, where
$c_{ij}$ is constant and ${\mathbb C}\ni z=x+\sqrt{-1}~y=(x,y)\in {\mathbb R}^2$.  In these charts, the
transverse measure
is given by integrating
$|dy|$ along transversals.

There is a natural equivalence relation on the set of all measured foliations in $F$, and there is a natural
topology on
the set of equivalence classes; see \cite{FLP} and \cite{PH}.  Roughly, if $c$ is an essential simple closed
curve in $F$
transverse to
${\cal F}$, one can
evaluate $\mu$ on $c$ to determine its ``geometric intersection number'' with ${\cal F}_\mu$; each homotopy class
$[c]$ of such curve has a representative minimizing this intersection number, and this minimum value is called
the  {\it
geometric intersection number} $i_{{\mathcal F}_\mu}[c]$ of $[c]$ and ${\cal F}_\mu$.  This describes a mapping
${\mathcal
F_\mu}\mapsto i_{{\mathcal F}_\mu}$ from the set of measured foliations to the function space
${\mathbb R}_{\geq 0}^{\mathcal
S(F)}$, where
${\mathcal S}(F)$ is the set of all homotopy classes of essential simple closed curves in $F$, and the function
space is given the weak topology.   The equivalence classes of measured foliations can be described as the fibers of
this
map,
and the topology as the weakest one so that each $i _{\cdot}[c]$ is continuous.
(It does not go unnoticed that this effectively ``quantizes the observables corresponding
to closed curves''.)
Both the
equivalence relation and the topology can be described more geometrically; see
\cite{FLP} and \cite{PH}.  In fact, the equivalence relation on measured foliations is generated by isotopy
and ``Whitehead moves'', which are moves on measured foliations dual to those depicted in Figure~13a.

There is thus a space ${\mathcal M\mathcal F}(F)$ of measured foliation classes on $F$ embedded in
${\mathbb R}_{\geq 0}^{\mathcal S(F)}$.
Each measured foliation (class) ${\cal F}_\mu$ determines an underlying {\it projective measured foliation (class)},
where
one projectivizes $\mu$ by the natural action of ${\mathbb R}_{>0}$ on measures and obtains the space of
{\it projective measured foliations} ${\mathcal P\mathcal F}(F)$ as the quotient of
${\mathcal M\mathcal F}(F)$ by this action.

\vskip .1in

\noindent{\bf Remark}~
It is necessary later to be a bit formal about the empty foliation in $F$, which we shall
denote by ${ 0}$ and identify with the zero functional in ${\mathbb R}_{\geq 0}^{\mathcal S(F)}$.  Let
${\mathcal M\mathcal F}^+(F)$ denote ${\mathcal M\mathcal F}(F)$ together with ${ 0}$ topologized so that a
neighborhood of ${ 0}$ is homeomorphic to the cone from ${ 0}$ over ${\mathcal P\mathcal F}(F)$.  In other words,
${\mathcal P\mathcal F}(F)$ is the projectivization of the ${\mathbb R}_{>0}$-space ${\mathcal M\mathcal F}^+(F)-\{{
0}\}\subseteq{\mathbb R}_{\geq 0}^{\mathcal S(F)}$.

\vskip .1in

Projective measured foliations were introduced by William Thurston in the 1970's as
a tool for studying the degeneration of geometric structures in dimensions two and three as well as for studying the
dynamics of homeomorphisms in two dimensions.

If
$[d]\in {\mathcal S}(F)$, then there is a corresponding functional $i _{[d]}\in{\mathbb R}_{\geq 0}^{\mathcal S(F)}$,
where $i _{[d]}[c]$ is the minimum number of times that representatives $c$ and $d$ intersect, counted
{\sl without sign},
the
``geometric intersection number''. In effect, the space ${\mathcal P\mathcal F}(F)$ forms a completion of
${\mathcal S}(F)$; more precisely, the projective classes of $\{i_{[d]}\cdot :[d]\in {\mathcal S}(F)\}$
are dense in the projectivization of ${\mathcal M\mathcal F}^+(F)-\{{0}\}$.

Furthermore, ${\mathcal P\mathcal F}(F)$ is
homeomorphic to a piecewise-linear sphere of dimension $6g-7$.  This sphere compactifies the usual
Teichm\"uller space of $F$ so as to produce a closed ball of dimension $6g-6$ upon which the usual mapping
class group
$MC(F)$
of
$F$  acts continuously.  For instance, one immediately obtains non-trivial results from the Lefshetz fixed
point theorem.

Thurston's
boundary does not descend in any tractable geometric sense to the quotient by $MC(F)$ since
the $MC(F)$-orbit of any non-separating curve is dense in ${\mathcal PF}(F)$.
The quotient ${\mathcal PF}(F)/MC(F)$ is thus
dramatically  non-Hausdorff.

\vskip .2in

\subsection*{{\rm B.2} Measured foliations and c/o structure}

\vskip .2in

Unlike the body of the paper, where we strived to include only those combinatorial aspects which are
manifest for physical interactions of strings, here we briefly describe a more speculative mathematical extension of
the
foregoing theory in the context of general measured foliations in a windowed surface
$F=F_g^s(\delta _1,\ldots ,\delta _r)$ with windows
$W$, set
$\sigma$ of punctures of cardinality $s\geq 0$, and set  $\delta$ of distinguished points on the boundary.
 Let $\beta$
denote a brane-labeling on $F$, and set $\delta ({\beta})=\{ g\in \delta : \beta (g)\neq \emptyset\}$.

A  measured $\beta$-foliation of $F$ is a measured foliation
${\mathcal F}_\mu$ in the usual sense of a closed subsurface (perhaps with boundary or punctures)  of $F$ so
that leaves
of ${\mathcal F}$ are either simple closed curves (which may be neither contractible nor
puncture-parallel nor boundary-parallel),  bi-infinite lines, or line segments with endpoints in
$\cup W$ which are not boundary-parallel in $F-\delta (\beta )$.

We shall furthermore require that ${\mathcal F}_\mu$ has
{\it compact
support} in the sense that its leaves are disjoint from a neighborhood of $\delta\cup\sigma$.  In particular,
${\mathcal F}$ is not permitted to have leaves that are asymptotic to $\delta\cup\sigma$.

There is again a natural equivalence relation on the set of all measured $\beta$-foliations of compact support
in $F$ and
a  natural topology on the space of equivalence classes induced by geometric intersection numbers with curves
as before
and now
also with embedded arcs connecting points of
$\delta(\beta )$.
Let $${\mathcal M\mathcal F}_0(F,\beta )$$ denote the corresponding space of measured
$\beta$-foliations of compact support.

We shall go a step further and allow puncture- and boundary-parallel curves (for instance, in order to capture that
part of
the foliation possibly discarded in the body of the paper in open self-gluing):  a {\it non-negative collar weight}
on a
windowed surface is a ${\mathbb R}_{\geq 0}$-function defined on $\sigma$, and one imagines an annulus of
width given by the collar weight foliated by puncture- or boundary-parallel curves.
Let $$\widetilde{\mathcal M\mathcal F}_0^{\geq 0}(F,\beta )\approx {\mathcal M\mathcal F}_0(F,\beta )\times{\mathbb
R}_{\geq 0}^s$$ denote the corresponding space of measured
$\beta$-foliations of compact support together with a non-negative collar weight.

It is nearly a tautology that any partial measured foliation of compact support with non-negative collar weight
decomposes
uniquely into a disjoint union  of its
``minimal'' sets, which are  one of the following:  a band of leaves parallel to an arc properly embedded $F$ with
endpoints in
$\cup W$; an annulus in $F$ foliated by curves parallel to the core of the annulus; a
foliation disjoint from the boundary  with
no closed leaves.

Finally, define
$$\widetilde{\mathcal M\mathcal F}_0^{\geq 0}(n,m,s)=\Biggl (\bigsqcup _F~\bigsqcup _\beta ~\widetilde
{\mathcal M\mathcal
F}_0^{\geq 0}(F)
\biggr
)/\sim,$$ where the outer disjoint union is over all homeomorphism classes of windowed surfaces $F$ with
$s\geq 0$ punctures, $n\geq 0$ closed windows, and $m\geq 0$ open windows with $m+n+s>0$,
and the inner disjoint union is over all  brane-labelings $\beta$ on $F$; the equivalence relation $\sim$ on the
double
disjoint union is generated by the  following identifications:  If ${\mathcal F}_\mu$ is a partial measured
foliation of
$F$ and
$\partial$ is a boundary component of
$F$ containing no active windows, then we collapse $\partial$ to a new puncture brane-labeled by the union of the
labels
on
$\partial$ to produce in the natural way a measured foliation ${\mathcal F}'_{\mu '}$ of another brane-labeled
windowed
surface
$F'$, and we identify
${\mathcal F}_\mu$ in $F$ with ${\mathcal F}'_{\mu '}$ in $F'$ in
$\widetilde{\mathcal M\mathcal F}_0^{\geq 0}(n,m,s)$ in
the natural way. In particular, each equivalence class has a representative measured foliation ${\mathcal F}_\mu$ in
some
well-defined topological type of windowed surface $F$,
where every  boundary component of $F$ has at least one active window for ${\mathcal F}_\mu$.

The gluing operations in Section~2 extend naturally to corresponding operations on the objects $\widetilde{\mathcal
M\mathcal
F}_0^{\geq 0}(n,m,s)$, where we assume that the boundary component containing each closed window is framed. In
each case,
we may also glue inactive windows to inactive windows in analogy to Figure~5.  Furthermore, rather than discard the
part
of
${\cal F}_3$ (in the notation of Section~2) that is not a band of arcs, we discard only those annuli foliated by null
homotopic simple closed curves as may arise from closed gluing or self-gluing.

We know of no physical interpretation for the minimal sets of a measured foliation other than bands of
arcs as in the body of the paper.  Furthermore, the quotient of
$\widetilde{\cal MF}_0^{\geq 0}(F,\beta )$ by $MC(F)$ is typically non-Hausdorff, yet contains the Hausdorff subspace
$\widetilde{Arc}'(F,\beta )$ studied in the body of the paper.  (The natural appearance of a non-Hausdorff space
as part of this theory does not go unnoticed.)

Geometrically, minimal sets that are not bands of arcs ``serve to mitigate other interactions'' for the simple
 reason that
a foliated band of arcs cannot cross an annulus which is foliated by circles parallel to the core of the annulus.

\vskip .2in

\centerline{\epsffile{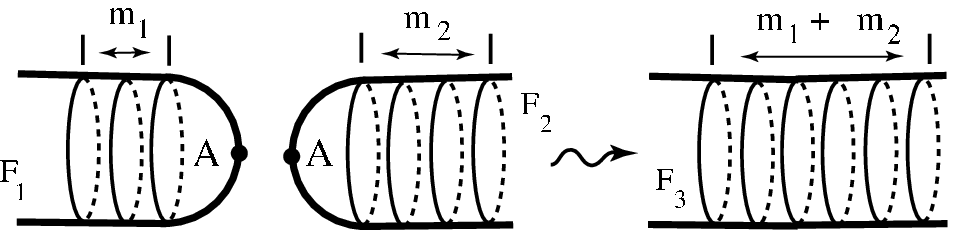}}

\vskip .2in

\centerline{\begin{footnotesize}{\bf Figure B.1}~Unpinching and self-unpinching.\end{footnotesize}}

\vskip .2in

There are also in the current context two further operations of ``unpinching'' and ``self-unpinching''  which are
geometrically natural and are illustrated in  Figure~B.1.

In the wider context of this appendix, minimal sets that are not
bands play three roles: they arise naturally from the  gluing
and self-gluing operations so as to mitigate other interactions; they arise from certain cases of open self-gluing
as
annuli
foliated by puncture-parallel curves, which themselves interact via unpinching; and they can be included as
{\sl a priori}
data which is invisible to the c/o structure on $\widetilde{Arc}(n,m)$ yet to which this c/o structure
contributes via
gluing and self-gluing.

Theorem~\ref{DT} holds essentially {\sl verbatim} in the current context (this was the
original purview) and describes the indecomposables as well as global coordinates.  There is, however, a more elegant
parametrization from
\cite{PP},  which should provide useful variables for quantization of the foregoing theory as closely related
coordinates did in Kashaev's quantization of decorated Teichm\"uller space
\cite{Ka}.  This parametrization arises by relaxing non-negativity of collar weights as follows.

Construct the space $\widetilde{\mathcal M\mathcal F}_0(F,\beta )$ in analogy to
$\widetilde{\mathcal M\mathcal F}_0^{\geq 0}(F,\beta )$ but allow the collar
weight on any puncture to be any real number; one imagines either a foliated annulus as before if the collar
weight is
non-negative or a kind of ``deficit'' foliated annulus if the collar weight is negative.

To explain the parametrization, let us fix a space-filling  brane-label $\beta
_A$ for simplicity, where $A\neq\emptyset$,  and choose a generalized pants
decomposition $\Pi$ of $F$. Inside each complementary region to
$\Pi$ insert a branched one-submanifold as illustrated in
Figure~B.2a, and combine these in the natural way to get a
branched one-submanifold $\tau$ properly embedded in $F$. Inside
each $F_{0,(3)}^0$ or $F_{0,1}^1$, there is a small triangle, as
illustrated, and the edges of these triangles are called the {\it
sectors} of $\tau$. Define a {\it measure} on $\tau$ to be the
assignment of a real number to each sector of $\tau$ subject to
the constraints that the ``coupling equations'' hold for each edge
of $\Pi$; namely, $a=b$ on $F_{0,1}^1$ and $a+b=c+d$ on
$F_{0,(4)}^0$ as illustrated in Figure~B.2b

\vskip .2in

\centerline{\epsffile{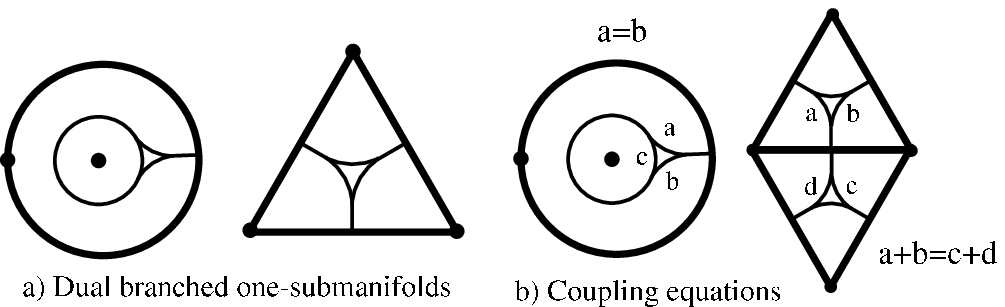}}

\vskip .2in

\centerline{\begin{footnotesize}{\bf Figure B.2} Dual branched one-submanifolds and coupling equations.
\end{footnotesize}}

\vskip .2in

One result from \cite{PP} is that the vector space of measures on $\tau$ is isomorphic to
$\widetilde{\mathcal M\mathcal F}_0(F,\beta _A)$, and there is a canonical fiber bundle
$$\widetilde{\mathcal M\mathcal F}_0(F,\beta _A )\to {\mathcal M\mathcal F}^+_0(F,\beta _A),$$
where the fiber over a point is given by the set of all collar weights on $F$,
and ${\mathcal M\mathcal F}^+_0(F,\beta _A)$ denotes the space of measured $\beta _A$-foliations of compact support
on $F$ completed by the empty foliation as in the Remark in the previous section.
Indeed, each puncture $p$
corresponds to a closed edge-path on $\tau$ which traverses a collection of sectors, and given a
measure $\mu$ on $\tau$, the collar weight $c_p$ of
$p$ is the minimum value that $\mu$ takes on these sectors; modify the original measure $\mu$ on $\tau$ by taking
$\mu '(b)=\mu (b)-c_p$ if $b$ is contained in the closed edge-path for $p$.  Since $\mu$ satisfies the coupling
equations, $\mu '$ extends uniquely to a well-defined non-negative measure on $\tau$, which describes a (possibly
empty)
element of ${\mathcal M\mathcal F}^+_0(F,\beta _A)$ by gluing together bands as before, one band for each edge of
$\tau$
on which
$\mu '>0$. A PL section of this bundle gives a PL embedding of the PL space ${\mathcal M\mathcal F}^+_0(F,\beta _A)$
into
the vector space $\widetilde{\mathcal M\mathcal F}_0(F,\beta _A )$.

Another result from
\cite{PP} for $r=0$ and $s\neq 0$ is that the Weil-Petersson K\"ahler two-form on Teichm\"uller space extends
continuously
to the natural symplectic structure given by Thurston on the space of measured foliations of compact support.


\vskip .3in

\renewcommand{\theequation}{C-\arabic{equation}}
\renewcommand{\thesection}{C}
\setcounter{equation}{0}  
\setcounter{subsection}{0}
\section*{Appendix C: Flows on $\widetilde{Arc}(n,m)$}

\vskip .2in

In this appendix, we shall define and study useful flows on $\widetilde{Arc}(n,m)$, two flows for
each window.  Fix a brane-labeled windowed surface $(F,\beta )$
with distinguished window $w$.  If $\alpha\in \widetilde{Arc}(F,\beta )$, then we shall now denote the $\alpha$-weight
of $w$ simply by
$\alpha (w)$.

The first flow $\psi ^w_t$ is relatively simple to define:
$$\psi _t^w(\alpha )=\bigl (1-t+{t/{\alpha (w)}}\bigr )\cdot \alpha ,$$
where the multiplication $x\cdot\alpha $ scales all the weights of $\alpha$ by the factor $x$;
thus, $\psi _0^w$ is the identity, and $\psi _1^w(\alpha ) (w)=1$.
This flow $\psi ^w_t$ provides a rough paradigm for the more complicated one to follow, and it alone is enough
to define open and closed {\sl non self-gluing} as follows:  given families $a$ in $\A(F_1,\beta _1)$
and $b$ in $\A(F_2,\beta _2)$, define
$$
a\odot _{v,w} b=
\begin{cases}
\psi ^v_1(a)\circ _{v,w}\psi ^w_1(b),&~{\rm if}~u,v~{\rm are~closed~windows;}\\
\psi ^v_1(a)\bullet _{v,w}\psi ^w_1(b),&~{\rm if}~u,v~{\rm are~open~windows.}
\end{cases}
$$
where $v,w$ are respective distinguished windows of $F_1,F_2$.
The gluing on the right is defined on the chain level provided $v,w$ are either both open or both closed
since $\psi _1^v(a)(v)\equiv 1\equiv \psi _1^w(b)(w)$.

\vskip .2in

\begin{lem}\label{clem1}
The operations $\odot _{v,w}$ on chains descend to well-defined
operations on homology classes.  Furthermore, suppose that $a,b,c$
are respective parameterized families in the brane-labeled
surfaces $(F_i,\beta _i )$, for $i=1,2,3$, with distinguished
windows $u$ in $F_1$, $v\neq v'$ in $F_2$, and $w$ in $F_3$, where
$\{ u,v\}$ and $\{ v',w\} $ each consists of either two open or
two closed windows.  Then there is a canonical homotopy between
$(a\odot _{u,v} b)\odot _{v',w} c$ and $a\odot _{u,v} (b\odot
_{v',w} c)$.
\end{lem}

\vskip .2in

\begin{proof}
To be explicit in this context of parameterized families, to say
that two families $a_0,a_1$ in $\A(F_1,\beta _1)$ of degree $k$
are homologous means that there is a degree $k+1$ family $A$ in
$\A(F_1,\beta _1)$ so that the boundary of the parameter domain
for $A$ decomposes into two sets $I_0,I_1$ with disjoint interiors
so that $A$ restricts to $a_i$ on $I_i$, for $i=0,1$.  It follows
that $\psi ^v_1 A\circ _{v,w} \psi ^w_1 b_0$ gives the required
homology between $a_0\odot _{v,w} b_0$ and $a_1\odot _{v,w} b_0$,
for any family $b_0$ in $\A(F_2,\beta _2)$.  The analogous
argument applies to two homologous families $b_0,b_1$ in
$\A(F_2,\beta _2)$, so $a_1\odot _{v,w} b_0$ is likewise
homologous to $a_1\odot _{v,w} b_1$. Thus, $a_0\odot _{v,w} b_0$
and $a_1\odot _{v,w} b_1$ are indeed homologous, completing the
proof that the operations are well-defined on homology.

As for the canonical homotopy, we claim that $(a\odot _{u,v} b)\odot _{v',w} c$ and
$a\odot _{u,v} (b\odot _{v',w} c)$ represent the same projective class.  Specifically,
let $F_{12}$ denote the surface containing $a\odot _{u,v} b$, let $F_{23}$ denote the surface
containing $b\odot _{v',w} c$, and let $F_{123}$ denote the common surface containing
$(a\odot _{u,v} b)\odot _{v',w} c$ and
$a\odot _{u,v} (b\odot _{v',w} c)$ with its induced brane-label $\beta _{123}$.  Corresponding to $\{ u,v\}$ in
$F_{123}$ there is either a properly embedded arc (if $u,v$ are open) or perhaps a simple closed curve (if $u,v$ are
closed or under certain circumstances if $u,v$ are open), and likewise corresponding to $\{ v',w\}$, there is an arc
or curve.   In the family $(a\odot _{u,v} b)\odot _{v',w} c$, the latter arc or curve has transverse measure constant
equal to one while the former arc or curve has some constant transverse measure $x$; in the family $a\odot _{u,v}
(b\odot _{v',w} c)$, the former arc or curve has transverse measure constant equal to
one while the latter arc or curve has some constant transverse measure $y$.
It follows from the definition of composition $\odot _{v,w}$ that $y=1/x$ and
$$x\cdot\bigl [ a\odot _{v,w}(b\odot _{v',w} c )\bigr ] = (a\odot _{v,w} b)\odot _{v',w} c,$$
where again $\cdot$ denotes the natural scaling action of ${\mathbb R}_{>0}$ on arc families in
$F_{123}$.  The required homotopy $\Psi _t$, for $0\leq t\leq 1,$ is finally given by
$$(1+t(x-1))\cdot\bigl [ a\odot_{v,w}(b\odot _{v',w} c)\bigr ].$$
\end{proof}

\vskip .2in

Turning now to preparations for the second more intricate flow for self-gluings,
suppose that $F$ is a windowed surface
with brane-labeling $\beta$.  If $\alpha\in Arc(F,\beta )$ and $w$ is some specified window of $F$, then the
bands of $\alpha$ that meet $w$ can be grouped together as follows: consecutive bands along $w$ that
connect
$w$ to a common window $w'$ are grouped together into the {\it $w$-bands} of $\alpha$ at $w$, where
consecutive
$w$-bands are not permitted to share a common endpoint other than $w$ in the obvious terminology.  In particular, for
any window
$w'$ of
$F$, there is the collection of $w$-bands of $\alpha$ with endpoints $w$ and $w'$.  Still more particularly, there
are the $w$-bands that have both endpoints at $w$, which are called the {\it self bands} of $w$.

Fix windows $w$ and $w'\neq w$ of $F$ and assume that $\alpha\in\A (F,\beta )$ satisfies $\alpha (w)>\alpha (w')$.  We
shall define a flow $\phi ^w_t$, for $-1\leq t< 1$,  so that at a certain first critical time $t_c<1$, we have
equality $\phi ^w_{t_c} (w)=\phi ^w_{t_c} (w')$.

The flow $\phi ^w_t$ is defined in two stages for $-1\leq t\leq 0$ and for $0\leq t<1$, and the first stage is
relatively easy to describe:  leave alone the $w$-bands of $\alpha$ other than the self bands, and scale
the weight of each self band by the factor $|t|$.  Thus, $\phi _{-1}^w$ is the identity, and
$\phi _{0}^w(\alpha )$ has no self bands at $w$.  Furthermore, $\phi ^w_t(\alpha) (w)$ is monotone decreasing, and
$\phi ^w_t(\alpha ) (u)$ is constant independent of $t$ for any $u\neq w$.

If $\alpha$ has only self bands at $w$, then $\phi _{0}^w(\alpha )(w)=0$ (so the flow is defined in
$\A(F,\beta )$ only for $-1\leq t<0$), and there is thus
some smallest
$t_c<0$ so that $\phi ^w_{t_c}(\alpha ) (w)=\phi ^w_{t_c} (\alpha )(w')$.  More generally, even if $\alpha$ has non
self bands at
$w$, it may happen that there is some smallest $t_c\leq 0$ so that $\phi ^w_{t_c} (\alpha )(w)=\phi ^w_{t_c}
(\alpha )(w')$.  This completes the definition of the first stage of the flow up to the point that there are no self
bands at $w$.

If there is no such
$t_c\leq 0$, then we continue to define the second stage of the flow $\phi ^w_t$ for $0\leq t<1$ in the absence of
self bands as follows.

If there is only one $w$-band of weight $a$, then $\phi ^w_t$ is defined as in Figure~C.1, where the
darkened
central part of the original foliated rectangle corresponding to this $w$-band is left alone, and the outer
white part of the foliation is erased, i.e., leaves are removed from the foliation.

\vskip .2in

\centerline{\epsffile{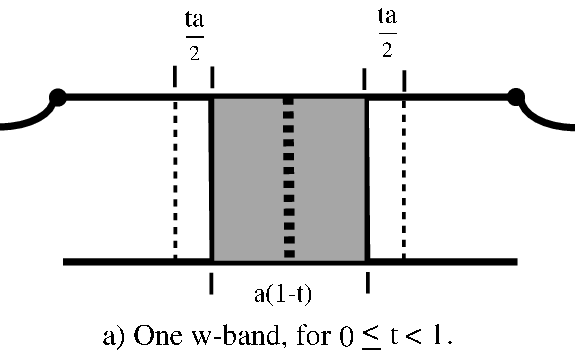}}

\begin{footnotesize}
\centerline {{\bf Figure~C.1} The flow $\phi^w_t$ for one w-band.}
\end{footnotesize}

\vskip .2in

\centerline{\epsffile{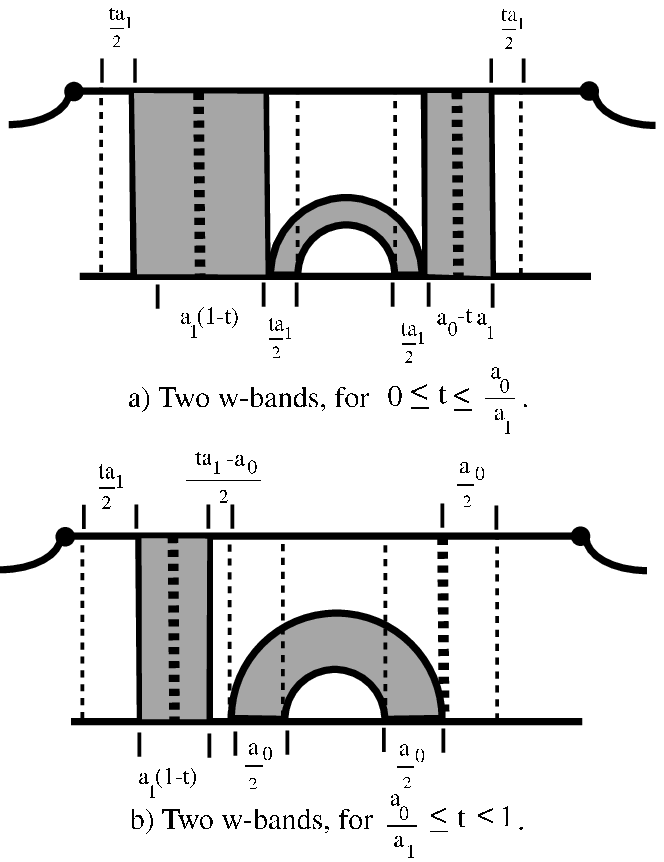}}

\begin{footnotesize}
\centerline {{\bf Figure~C.2} The flow $\phi^w_t$ for two w-bands.}
\end{footnotesize}

\vskip .2in

More
interesting is the case that there are two $w$-bands, which is illustrated in Figure~C.2.  In between the two bands,
we
surger together arcs preserving measure in the manner indicated.  Since the $w$-bands are consecutive and there are
no self bands at $w$, the resulting  arcs
must connect distinct windows, hence must be essential and moreover cannot be a self band at any window.  We erase
leaves from the other sides of the two bands as before. Letting $a_0<a_1$ denote the weights of the two $w$-bands,
there is a critical time
$t=a_0/a_1$ when there  is
a unique $w$-band.  The flow before the critical time is illustrated in Figure~C.2a and after it in Figure~C.2b.

\vskip .2in

\centerline{\epsffile{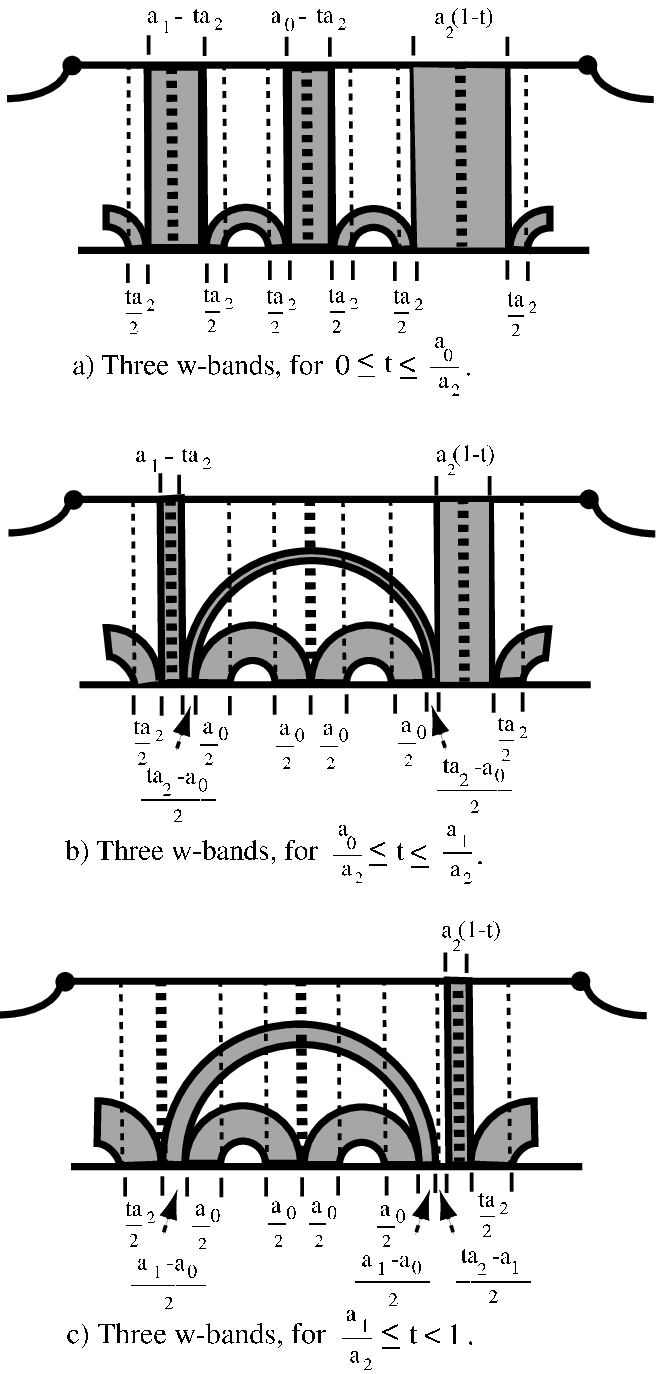}}

\begin{footnotesize}
\centerline {{\bf Figure~C.3} The flow $\phi^w_t$ for three or more w-bands, first case.}
\end{footnotesize}

\vskip .2in

For three or more bands, there are two essential cases.  Let us fix three consecutive $w$-bands $b_1,b_2,b_3$
of respective weights $a_1,a_0,a_2$, where $a_0\leq a_1\leq a_2$, i.e., $b_2$ is of minimum weight and $b_1\leq
b_3$.  It may happen  that
$b_1$ and $b_3$ do not share an endpoint other than $w$, and in this case, the flow is defined as illustrated in
Figure~C.3.  There are two critical times $t=a_0/a_1, a_1/a_2$, at each of which the number of $w$-bands is
decreased by one.

Figure~C.4 illustrates the case that $b_1$ and $b_3$ do share an endpoint other than $w$.  At the first
critical time
$t=a_0/a_1$, the number of bands is in effect decreased by two since the two bands must now be combined to one; at
the
second critial time $t=a_1/a_2$, there is then a transition from two-bands-as-one to a single band.

\vskip  .2in

\centerline{\epsffile{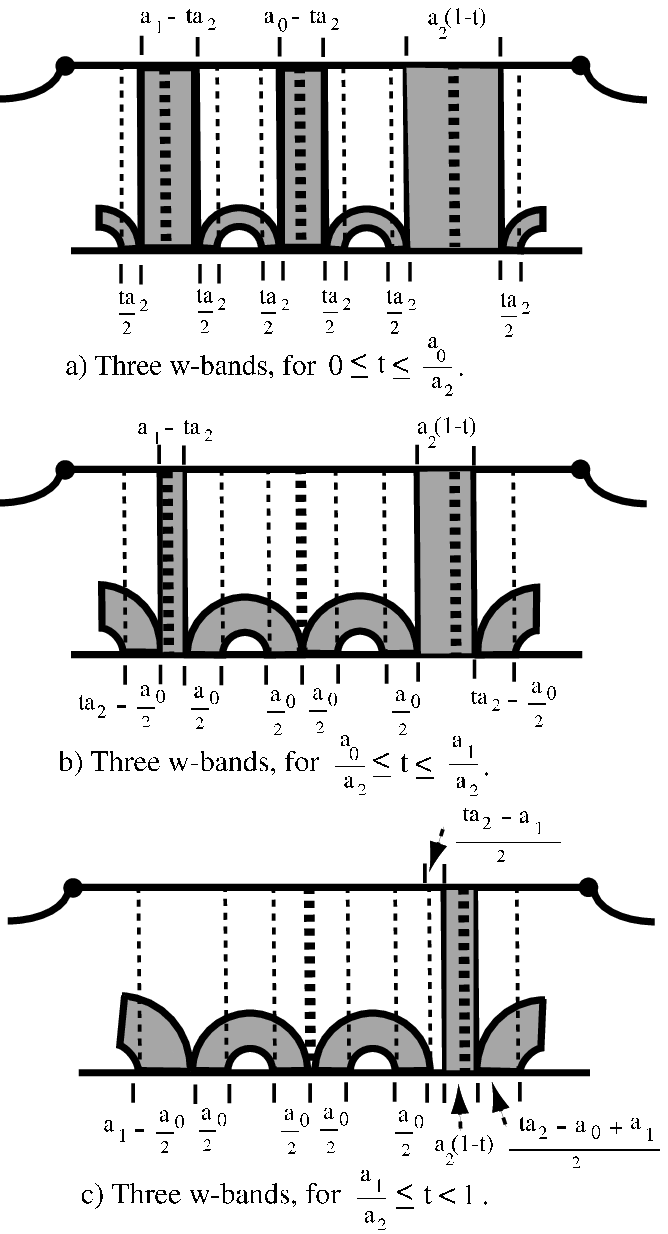}}

\begin{footnotesize}
\centerline {{\bf Figure~C.4} The flow $\phi^w_t$ for three or more w-bands, second case.}
\end{footnotesize}

\vskip .2in

Now generalizing Figures~C.3 and C.4 in the natural way to consecutive $w$-bands, this completes the definition of
the
flow
$\phi ^w_t$.  Notice that additive relations among the weights of the $w$-bands can lead to modifications of the
evolution, but in all cases, there are two basic types of critical times when a band becomes exhausted: either the
newly consecutive
$w$-bands share an endpoint other than $w$ so must be combined to a single $w$-band, or they do not and one adds a
new
band of surgered arcs which does not meet $w$.

By definition, $\phi _0^w$ is the identity.  Furthermore for any $\alpha\in\A(F,\beta )$
by construction, we have $\phi _t^w(\alpha ) (w)$ tending to zero as $t$ tends to one.
In fact, consideration of the formulas in Figures~C.1-C.4 shows that
${d/{dt}}~\phi _t^w(\alpha )(w)= -\phi _t^w(\alpha )(w)$, so the decay to zero
is exponential in $t$.

\vskip .2in

\begin{lem}\label{clem2}Fix a brane-labeling $\beta$ on a windowed surface $F$ with respective closed and open windows
$S$ and
$T$, fix a window
$w$ of
$F$, and let
$\A_w(F,\beta)\subseteq \A(F,\beta )$ denote the subspace corresponding to arc families $\alpha$
that have no self bands at $w$.
Then there is a
continuous piecewise-linear flow
$$\phi _t^w:\A_w(F,\beta )\to \A_w(F,\beta ),~{\rm for}~0\leq t<1,$$
such that $\phi_0^w$ is the identity, and for any other $w\neq w'\in S\sqcup T$ with
$\alpha (w')\leq \alpha (w)$, there is a first time $t_c=t_c(\alpha)<1$ for which
$$\phi ^w_{t_c}(\alpha ) (w)=\phi ^w_{t_c}(\alpha ) (w'),$$
where $t_c(\alpha )$ depends continuously on $\alpha$.  \end{lem}

\vskip .2in

\begin{proof}
Suppose that $w\neq w'$ is another window of $F$, and consider the $w'$-bands of $\alpha$.  If
there is a $w'$-band with an endpoint distinct from $w$, then $\phi _t^w(\alpha )(w')$ is bounded below
uniformly in $t$ since such a $w'$-band is undisturbed by the flow $\phi _t^w$ by definition.
Furthermore, suppose that two consecutive $w$-bands have respective other endpoints $w'$ and $w''\neq
w$.  After an instant of time, the flow combines these bands to produce a $w'$-band as in the previous
sentence.  It follows that if $\a(w')\leq \a (w)$, then there is indeed some critical first time
$t_c=t_c(\alpha )=t_c(\alpha, w')$ so that
$\phi_{t_c}(\a)(w)=\phi_{t_c}(\a)(w')$.

As for continuity, a neighborhood in $\widetilde{Arc}(F,\beta )$
of a weighting $\mu$ on an arc family $\alpha$ is a choice of
maximal arc family $\alpha '\supseteq\alpha$ together with an open
set $V$ of weightings $\mu '$ on $\alpha '$ so that $\mu '\in V$
restricts to a neighborhood of $\mu $, and values of $\mu '$ on
$\alpha '-\alpha$ are all bounded above by some $\varepsilon$.
We may choose a maximal arc family $\alpha '$ with no self bands at $w$ since
$\alpha$ has no self bands at $w$.
Adding arcs in $\alpha '-\alpha$ to $\alpha$ cannot decrease the
number of $w$-bands.  Furthermore, if $a_{max}$ is the maximum
value of $\alpha$ on a $w$-band, then in time $t>\varepsilon
~a_{max} $, the corresponding bands of arcs in $\alpha '$ are
erased or combined by $\phi ^w_t(\alpha ')$. Continuity of the flow follows from
these facts.

Continuity of the critical time function also follows from these
facts and the following considerations. The number of $w$-bands
for $\phi ^w_t(\alpha )$ is a non-increasing function  of $t$, so
either $\phi _{t_0}^w(\alpha )(w)=\phi _{t_0}^w(\alpha )(w')$ for
some $t_0<t$, or there is a unique $w$-band at time $t$.  In the
latter case, either this $w$-band has its other endpoint distinct
from $w'$, so $\phi _{t}^w(\alpha )(w)$ is exponentially
decreasing in $t$ while $\phi _{t}^w(\alpha ) (w')$ is constant,
or this $w$-band has its other endpoint at $w'$.  In the latter
case, either this $w$-band is not the unique $w'$-band, so $\phi
_t^w(w')>\phi _t^w(\alpha )(w)$ and these quantities must have
therefore agreed at some time preceding $t$, or there is a unique
$w$-band and a unique $w'$-band connecting $w$ and $w'$.  In the
latter case, $\phi _t^w(\alpha )(w')=\phi _t^w(\alpha )(w)$ for
all time after $t$.
\end{proof}

\vskip .2in

As before, suppose that $a=a(s)$ is a parameterized family in
$\A(F_1,\beta _1)$ and $b=b(t)$ is a parameterized family in
$\A(F_2,\beta _2)$, where $v,w$ are respective distinguished
windows of $F_1, F_2$.  Let

\begin{eqnarray*}
u(s,t)&=&
\begin{cases}
v,&~{\rm if}~a(s)(v)\geq b(t)(w);\\
w,&~{\rm if}~a(s)(v)< b(t)(w);
\end{cases}\\
d(s,t)&=&
\begin{cases}
a(s),&~{\rm if}~a(s)(v)\geq b(t)(w);\\
b(t),&~{\rm if}~a(s)(v)< b(t)(w),
\end{cases}
\end{eqnarray*}

\noindent and finally define

\begin{eqnarray*}\label{ceqn2}
a\odot ^{v,w} b=
\begin{cases}
\phi ^{u(s,t)}_{t_c} a(s)\circ ^{v,w}\phi ^{u(s,t)}_{t_c} b(t),&{\rm if}~v,w~{\rm
are~closed~windows;}\\
\phi ^{u(s,t)}_{t_c} a(s)\bullet ^{v,w}\phi ^{u(s,t)}_{t_c} b(t),&~{\rm if}~v,w~{\rm
are~open~windows,}
\end{cases}
\end{eqnarray*}

\noindent where $t_c=t_c(d(s,t))$ is the critical first time when the flow $\phi ^{u(s,t)}_{t_c}$ achieves
$a(s)(v)=b(t)(w)$ so that gluing is possible.

\vskip .2in

\begin{lem}\label{clem3}The operations $\odot ^{v,w}$ on chains descend to well-defined operations on homology
classes. Furthermore, suppose that $a,b,c$ are respective
parameterized families in the brane-labeled surfaces $(F_i,\beta
_i )$, for $i=1,2,3$, with distinguished windows $u$ in $F_1$,
$v\neq v'$ in $F_2$, and $w$ in $F_3$, where $\{ u,v\}$ and $\{
v',w\} $ each consists of either two open or two closed windows.
Then there is a canonical homotopy between $(a\odot ^{u,v} b)\odot
^{v',w} c$ and $a\odot ^{u,v} (b\odot ^{v',w} c)$.
\end{lem}

\vskip .2in

\begin{proof}
That the operations are well-defined on the level of homology follows in analogy to the previous case
Lemma~\ref{clem1}, and it remains only to describe the canonical homotopies.  To this end in addition to the
surfaces $F_{12}$ and $F_{23}$ respectively containing
$a\odot ^{u,v} b$ and $b\odot ^{v',w} c$, as well as the surface $F_{123}$ containing both
$(a\odot ^{u,v} b)\odot ^{v',w} c$ and
$a\odot ^{u,v} (b\odot ^{v',w} c)$, we must introduce another
auxiliary surface $F$ defined as follows.  Among the two operations $\odot ^{u,v}$ and $\odot ^{v',w}$, suppose
that $\kappa =0,1,2$ of the operations are closed string self gluings.  The auxiliary surface $F$ is homeomorphic to
$F_{123}$ except that it has $2\kappa$ additional punctures, which we imagine as lying in a small annular neighborhood
of the corresponding curve in $F_{123}$ with one new puncture on each side of the curve.

In this surface $F$ with punctures $\sigma$ and distinguished points $\delta$ on the boundary, we shall consider
collections of arc families somewhat more general than before.   Specifically, we shall now allow arcs to
have one or both of their endpoints in $\delta\cup\sigma$.  Arc families are still defined rel $\delta\cup\sigma$ as
before, and the geometric realization of the corresponding partially ordered set is the space within which we shall
define the required homotopy.  Arcs with endpoints at $\delta\cup\sigma$ are called ``special'' arcs, and there is a
homotopy that simply scales their weights to zero to produce an arc family in the usual sense in $F_{123}$; our
homotopies will be described in this augmented arc complex of $F$ taking care to make sure that this projection to
$F_{123}$ lies in $\A(F_{123},\beta _{123})$, i.e., every window has positive weight for the projection, where $\beta
_{123}$ is the brane-labeling induced on $F_{123}$ from the given data.

We modify the constructions of Figure~C.1-4 in one manner for open string self-gluing and in another manner for closed
string self-gluing.  For the former in the second stage of the homotopy, instead of erasing the outermost
edges of the one or two outermost $w$-bands, let us instead keep these arcs and run them as special arcs to the nearby
points of
$\delta$ in the natural way. For closed gluing, we employ the additional punctures of $F$ and instead of erasing the
outermost edges of the one or two outermost $w$-bands, we instead run them as special arcs to the nearby additional
punctures in the natural way.

There is thus a modified flow for the augmented arc families with the advantage that only the weight of
the window $w$ changes under the modified flow,  and therefore the modified flows corresponding to different windows
commute.  There is furthermore a modified operation defined in analogy
to $\odot ^{v,w}$ using the modified flow.  Because the modified flows of different windows commute, the modified
expressions
$(a\odot ^{u,v} b)\odot ^{v',w} c$ and $a\odot ^{u,v} (b\odot ^{v',w} c)$ agree exactly.

Finally, there are the special bands that arise from
the operation $\odot ^{u,v}$, and then there are the special bands that arise from the operation
$\odot ^{v',w}$.  Scaling the former to zero projects to one order of composition, and scaling the latter to zero
projects to the other order of composition.  This establishes the asserted homotopies.
\end{proof}

\vskip .2in

\begin{cor}\label{integrable}
For any brane-labeled windowed surface $(F,\beta )$, the space $\A(F,\beta )$ supports
the collection $\{ \phi ^w_t:w~{\rm is~a~window~of}~F\}$ of pairwise commuting flows.
\end{cor}

\vskip .2in

\begin{proof}
As in the previous proof, scaling to zero first the special bands of one modified flow and then scaling to zero the
special bands of the other modified flow projects to one order of composition of flows on $\A(F,\beta )$, while
scaling to zero in the other order produces the other order of composition.
\end{proof}

\vskip .2in

\begin{lem}\label{clem4}
Under the hypotheses and notation of Lemma~\ref{clem3}, there is a canonical homotopy between
$(a\odot ^{u,v} b)\odot _{v',w} c$ and
$a\odot ^{u,v} (b\odot _{v',w} c)$ and a canonical homotopy between
$(a\odot _{u,v} b)\odot ^{v',w} c$ and
$a\odot _{u,v} (b\odot ^{v',w} c)$
\end{lem}

\vskip .2in

\begin{proof}
As in Lemma~\ref{clem3}, there is an auxiliary surface with arc
families augmented by special arcs. Again, modified flows give
rise to modified operations so that the asserted pairs are
projectively equivalent.  Finally, scaling with homotopies as in
Lemma~\ref{clem1} completes the proof.
\end{proof}

\section*{Acknowledgments}
Both authors are happy to acknowledge the kind support of the Max- Planck Institute for Mathematics (Bonn) during
the Summer of 2005, where and when this work began, and particularly the warm hospitality of Yuri and Xenia Manin.
It is
also a pleasure for RCP to thank his colleagues Clifford Johnson and Tobias Ekholm for especially helpful discussions.

\end{document}